\documentclass[11pt,a4paper]{article}
\usepackage{amsfonts,amsmath,amssymb,accents,bm,amsthm}
\usepackage{verbatim}

\usepackage{hyperref}
\usepackage[normalem]{ulem}

\usepackage{bbm}

\usepackage{diagbox}

\usepackage{graphicx}
\usepackage{xcolor}

\newcommand{\dis}{\displaystyle}

\newcommand{\bigboxdot}{
  \mathop{
    \vphantom{\bigoplus} 
    \mathchoice
      {\vcenter{\hbox{\resizebox{\widthof{$\displaystyle\bigoplus$}}{!}{$\boxdot$}}}}
      {\vcenter{\hbox{\resizebox{\widthof{$\bigoplus$}}{!}{$\boxdot$}}}}
      {\vcenter{\hbox{\resizebox{\widthof{$\scriptstyle\oplus$}}{!}{$\boxdot$}}}}
      {\vcenter{\hbox{\resizebox{\widthof{$\scriptscriptstyle\oplus$}}{!}{$\boxdot$}}}}
  }\displaylimits 
}

\newcommand{\noi}{\noindent}
\newcommand{\halmos}{\rule{1ex}{1.4ex}}
\newcommand{\QED}{\nopagebreak{\hspace*{\fill}$\halmos$\medskip}}

\newcommand{\quand}{\quad\mbox{and}\quad}

\newtheoremstyle{mythm}
  {}
  {}
  {\itshape}
  {}
  {\bfseries}
  {}
  {.5em}
  {#1 #2 \thmnote{(#3)}}

\theoremstyle{mythm}
\newtheorem{theorem}{Theorem}
\newtheorem{proposition}[theorem]{Proposition}
\newtheorem{lemma}[theorem]{Lemma}

\newtheorem{exercise}[theorem]{Exercise}
\newtheorem{corollary}[theorem]{Corollary}
\newtheorem{conjecture}[theorem]{Conjecture}

\newtheorem{counterex}[theorem]{Counterexample}

\newcommand{\bt}{\begin{theorem}}
\newcommand{\et}{\end{theorem}}
\newcommand{\bl}{\begin{lemma}}
\newcommand{\el}{\end{lemma}}
\newcommand{\bp}{\begin{proposition}}
\newcommand{\ep}{\end{proposition}}
\newcommand{\bcor}{\begin{corollary}}
\newcommand{\ecor}{\end{corollary}}
\newcommand{\br}{\begin{remark}\rm}
\newcommand{\er}{\end{remark}}
\newcommand{\bcon}{\begin{conjecture}}
\newcommand{\econ}{\end{conjecture}}
\newcommand{\bex}{\begin{exercise}}
\newcommand{\eex}{\end{exercise}}
\newcommand{\bcou}{\begin{counterex}}
\newcommand{\ecou}{\end{counterex}}

\theoremstyle{definition}

\newtheorem{remark}[theorem]{Remark}

\newenvironment{Proof}[1][]{\noi\textbf{Proof #1}}{\QED}
\newcommand{\bpro}{\begin{Proof}}
\newcommand{\epro}{\end{Proof}}

\newcommand{\be}{\begin{equation}}
\newcommand{\ee}{\end{equation}}
\newcommand{\ba}{\begin{array}}
\newcommand{\ea}{\end{array}}
\newcommand{\bac}{\begin{array}{r@{\,}c@{\,}l}}
\newcommand{\bc}{\be\begin{array}{r@{\,}c@{\,}l}}
\newcommand{\ec}{\end{array}\ee}

\newcommand{\ga}{\gamma}

\newcommand{\de}{\delta}
\newcommand{\De}{\Delta}
\newcommand{\eps}{\varepsilon}
\newcommand{\la}{\lambda}
\newcommand{\La}{\Lambda}

\newcommand{\tet}{\theta}

\newcommand{\Ai}{{\cal A}}
\newcommand{\Bi}{{\cal B}}

\newcommand{\Hi}{{\cal H}}

\newcommand{\Ni}{{\cal N}}

\newcommand{\Ti}{{\cal T}}

\newcommand{\A}{{\mathbb A}}

\newcommand{\E}{{\mathbb E}}

\renewcommand{\P}{{\mathbb P}}

\newcommand{\R}{{\mathbb R}}

\newcommand{\Z}{{\mathbb Z}}

\newcommand{\sub}{\subset}

\newcommand{\Asto}[1]{\underset{{#1}\to\infty}{\Longrightarrow}}

\newcommand{\dgg}{\dagger}
\newcommand{\ov}{\overline}
\newcommand{\un}{\underline}

\newcommand{\ffrac}[2]{{\textstyle\frac{{#1}}{{#2}}}}

\newcommand{\di}{\mathrm{d}}

\newcommand{\ha}{\ffrac{1}{2}}

\newcommand{\psib}{{\boldsymbol\psi}}

\begin{document}

\makeatletter\@addtoreset{equation}{section}
\makeatother\def\theequation{\thesection.\arabic{equation}}

\renewcommand{\labelenumi}{{\rm (\roman{enumi})}}
\renewcommand{\theenumi}{\roman{enumi}}

\title{\vspace*{-2cm}Applying monoid duality to a double contact process}
\author{Jan Niklas Latz\footnote{The Czech Academy of Sciences,
  Institute of Information Theory and Automation,
  Pod vod\'arenskou v\v e\v z\' i~4,
  18200 Praha 8,
  Czech Republic.
  latz@utia.cas.cz, swart@utia.cas.cz.}
\and
Jan~M.~Swart${}^\ast$}

\date{\today}

\maketitle

\begin{abstract}\noindent
In this paper we use duality techniques to study a combination of the well-known \emph{contact process} (CP) and the somewhat less-known \emph{annihilating branching process}. As the latter can be seen as a \emph{cancellative} version of the contact process, we rebrand it as the \emph{cancellative contact process} (cCP). Our process of interest will consist of two entries, the first being a CP and the second being a cCP. We call this process the \emph{double contact process} (2CP) and prove that it has (depending on the model parameters) at most one invariant law under which ones are present in both processes. In particular, we can choose the model parameter in such a way that CP and cCP are monotonely coupled. In this case also the above mentioned invariant law will have the property that, under it, ones in the cCP can only be present at sites where there are also ones in the CP. Along the way we extend the dualities for Markov processes discovered in our paper ``Commutative monoid duality" to processes on infinite state spaces so that they, in particular, can be used for interacting particle systems.
\end{abstract}
\vspace{.5cm}

{\it Acknowledgements:} Work supported by grant 20-08468S of the Czech Science Foundation (GA\v{C}R).

{\setlength{\parskip}{-2pt}\tableofcontents}

\newpage

\section{Introduction}
\label{S:Int}

\subsection{Aim of the paper}

After having identified in \cite{LS22} a class of duality functions based on commutative monoids, our aim for this present paper is to apply one of those dualities to a concrete process. To do so we combine the contact process with its cancellative version, the process formerly known as the annihilating branching process. The considerations in \cite{LS22} indicate that this combined process has a self-duality that we use here to characterise all invariant laws of the process. 

To use the dualities discovered in \cite{LS22}, we first have to generalise the techniques presented in \cite{LS22} to infinite state spaces. This is done in Section~\ref{S:MProc} and is one of the main contributions of the present paper.

Additionally, in Section~\ref{S:Info}, we give precise definitions and some first results towards the goal of characterising all duality functions of the type considered in \cite{LS22} that determine the law of a process uniquely. This was posed as an open problem in \cite[Section~1.5]{LS22}.

\subsection{The processes of interest}
\label{S:PoI}

We set $T:=\{0,1\}$ and let $\Ti$ denote the space of all functions $x:\Z^d\to T$. Moreover, we let $\vee$ and $\oplus$ denote the binary operators on $T$ defined by the addition tables:

\begin{center}
\begin{tabular}{c|ccc}
$\vee$&0&1\\
\hline
0&0&1\\
1&1&1\\
\end{tabular}
\hspace{2cm}
\begin{tabular}{c|ccc}
$\oplus$&0&1\\
\hline
0&0&1\\
1&1&0\\
\end{tabular}
\end{center}

\noi
In words, this says that $x\vee y$ is the maximum of $x$ and $y$ and $x\oplus y$ is the sum of $x$ and $y$ modulo 2. For all $i,j\in\Z^d$, we define ``infection maps'' ${\tt inf}^\ast_{ij}:\Ti\to\Ti$ $(\ast\in\{\vee,\oplus\})$ and a ``death map'' ${\tt dth}_i:\Ti\to\Ti$ as follows:
\begin{align}
\label{genmaps}
\mathtt{inf}^\ast_{ij}(x)(k):=\begin{cases}
x(i)\ast x(j)&\text{if}\ k=j,\\
x(k)&\text{else},
\end{cases},\qquad
\mathtt{dth}_i(x)(k):=\begin{cases}
0&\text{if}\ k=i,\\
x(k)&\text{else}.
\end{cases}
\end{align} 
\noi 
We let $\Ni_i:=\{j\in\Z^d:\Vert i-j\Vert_1=1\}$ denote the set of nearest neighbours of a site $i\in\Z^d$ and write $i\sim j$ if $i$ and $j$ are nearest neighbours. We define formal generators
\begin{align}
\label{gener}
\dis G_{\ast}f(x):=\dis \lambda\sum_{i\in\Z^d}\sum_{j\in\Ni_i}
\big\{f\big({\tt inf}^\ast_{ij}(x)\big)-f\big(x\big)\big\}
+\delta\sum_{i\in\Z^d}\big\{f\big({\tt dth}_i(x)\big)-f\big(x\big)\big\}
\end{align}
for $\ast\in\{\vee,\oplus\}$, where $\lambda,\delta\geq0$ are model parameters. In words, we can describe the dynamics of the process generated by $G_\ast$ $(\ast\in\{\vee,\oplus\})$ as follows:
\begin{itemize}
\item At each site $i\in\Z^d$ sit two ``exponential clocks", one with rate $2d\lambda$ for \emph{reproduction} and one with rate $\delta$ for \emph{death}.
\item If the clock for reproduction at site $i\in\mathbb{Z}$ rings, the corresponding individual \emph{reproduces} by choosing a neighbouring site $j$ uniformly at random and adding its local state to the local state at $j$, where addition has to be interpreted in the sense of the operator $\ast$.
\item If the ``death clock" at site $i$ rings, individual $i$ \emph{dies} which means that its local state is replaced by 0, regardless of its previous value.
\end{itemize}
The process $C=(C_t)_{t\geq 0}$ with generator $G_\vee$ is the well-known contact process on $\Z^d$ with infection rate $\lambda$ and death rate $\delta$ (we denote this process shortly as CP($\lambda,\delta$)). The process $D=(D_t)_{t\geq 0}$ with generator $G_\oplus$, originally known as the annihilating branching process, we refer to as the cancellative contact process (cCP($\lambda,\delta$)). We chose the new name to stress the similarity of the processes, which differ only in the type of operator used in the definition of the infection maps ${\tt inf}^\ast_{ij}$ $(\ast\in\{\vee,\oplus\})$.

We will be interested in joint processes, consisting of a CP and a cCP, that are coupled in such a way that some of the infections and deaths happen for both processes at the same times. It will be helpful to write the generator of the coupled process in a form similar to (\ref{gener}). To achieve this, formally we define $U:=T\times T=\{0,1\}\times\{0,1\}$ and equip $U$ with $\veebar$, the product operator of $\vee$ and $\oplus$ from above, i.e.\ $(x,y)\veebar(v,w):=(x\vee v,y\oplus w)$ for $(x,y),(v,w)\in U$. This gives the following addition table:

\begin{center}
\begin{tabular}{c|cccc}
$\veebar$&(0,0)&(0,1)&(1,0)&(1,1)\\
\hline
(0,0)\rule[-0em]{0mm}{1.2em}&(0,0)&(0,1)&(1,0)&(1,1)\\
(0,1)\rule[-0em]{0mm}{1em}&(0,1)&(0,0)&(1,1)&(1,0)\\
(1,0)\rule[-0em]{0mm}{1em}&(1,0)&(1,1)&(1,0)&(1,1)\\
(1,1)\rule[-0em]{0mm}{1em}&(1,1)&(1,0)&(1,1)&(1,0)\\
\end{tabular}
\end{center}

\noi
In parallel to the above we denote by $\mathcal{U}$ the space of all functions $x=(x_1,x_2):\mathbb{Z}^d\to U$ and for each $i,j\in\Z^d$, we define infection maps ${\tt INF}_{ij},{\tt inf^1}_{ij},{\tt inf^2}_{ij}:\mathcal{U}\to\mathcal{U}$ and death maps ${\tt DTH}_i,{\tt dth^1}_i,{\tt dth^2}_i:\mathcal{U}\to\mathcal{U}$ as
\begin{align}
\begin{split}
\label{prodmaps}
&\mathtt{INF}_{ij}(x):=(\mathtt{inf}^\vee_{ij}(x_1),\mathtt{inf}_{ij}^\oplus(x_2)),\quad\mathtt{DTH}_i(x):=(\mathtt{dth}_i(x_1),\mathtt{dth}_i(x_2)),\\
&\mathtt{inf^1}_{ij}(x):=(\mathtt{inf}^\vee_{ij}(x_1),x_2),\quad\hspace{2.6em}\mathtt{dth^1}_i(x):=(\mathtt{dth}_i(x_1),x_2),\\
&\mathtt{inf^2}_{ij}(x):=(x_1,\mathtt{inf}^\oplus_{ij}(x_2)),\quad\hspace{2.6em}\mathtt{dth^2}_i(x):=(x_1,\mathtt{dth}_i(x_2)),\qquad(x=(x_1,x_2)\in\mathcal{U}),
\end{split}
\end{align} 
where the maps on the right hand sides are the maps from (\ref{genmaps}). We then define the generator $G_\veebar$ as
\begin{align}
\begin{split}
\label{2G}
G_\veebar f(x)&:=\lambda\sum_{i\in\Z^d}\sum_{j\in\Ni_i}
\big\{f\big({\tt INF}_{ij}(x)\big)-f\big(x\big)\big\}
+\delta\sum_{i\in\Z^d}\big\{f\big({\tt DTH}_i(x)\big)-f\big(x\big)\big\}\\
&\hspace{1.5em}+\lambda_\vee\sum_{i\in\Z^d}\sum_{j\in\Ni_i}
\big\{f\big({\tt inf^1}_{ij}(x)\big)-f\big(x\big)\big\}+\delta_\vee\sum_{i\in\Z^d}\big\{f\big({\tt dth^1}_i(x)\big)-f\big(x\big)\big\}\\
&\hspace{1.5em}+\lambda_\oplus\sum_{i\in\Z^d}\sum_{j\in\Ni_i}
\big\{f\big({\tt inf^2}_{ij}(x)\big)-f\big(x\big)\big\}+\delta_\oplus\sum_{i\in\Z^d}\big\{f\big({\tt dth^2}_i(x)\big)-f\big(x\big)\big\},
\end{split}
\end{align}
where $\lambda,\delta,\lambda_\vee,\delta_\vee,\lambda_\oplus,\delta_\oplus\geq0$ are model parameters. Standard results \cite[Theorem~4.30]{Swa22} tell us that the process $X=(X^1,X^2)=(X^1_t,X^2_t)_{t\geq0}$ with generator $G_\veebar$ is (like $C$ and $D$ before) well-defined. For later use, letting
\begin{align}
\label{S_fin}
\mathcal{U}_{\rm fin}:=\big\{x=(x_1,x_2)\in\mathcal{U}:\vert\{i\in\Z^d:(x_1(i),x_2(i))\neq(0,0)\}\vert<\infty\big\}
\end{align}
denote the set of finite configurations, one has, by Theorem~\ref{T:JansBIGtheo} below, for all choices of model parameters that
\be\ba{llll}\label{finpres}
\dis X_0\in\mathcal{U}_{\rm fin}\quad&\mbox{implies}\quad&\dis X_t\in\mathcal{U}_{\rm fin}\quad&\dis(t\geq 0)\quad{\rm almost\ surely.}
\ec
We call $X$ the \emph{double contact process} and denote it shortly as 2CP($\lambda,\delta,\lambda_\vee,\delta_\vee,\lambda_\oplus,\delta_\oplus$). If $X$ is a 2CP($\lambda,\delta,\lambda_\vee,\delta_\vee,\lambda_\oplus,\delta_\oplus$), then $X^1$ is a CP($\lambda+\lambda_\vee,\delta+\delta_\vee$) and $X^2$ is a cCP($\lambda+\lambda_\oplus,\delta+\delta_\oplus$). 

In particular, if $\lambda=\delta=0$, then $X^1$ and $X^2$ are independent processes. On the other extreme, if $\delta_\vee=\lambda_\vee=\delta_\oplus=\lambda_\oplus=0$, then $X^1$ and $X^2$ are fully coordinated in the sense that their infections and deaths happen at the same times. An interesting consequence of this choice of parameters is that the CP stochastically dominates the cCP. The first part of the following lemma says that this holds a bit more generally: if $\delta_\vee=\lambda_\oplus=0$ and the process is started in an initial state such that the CP dominates the cCP, then it follows from the definition of the maps in (\ref{prodmaps}) that this order is preserved by the evolution. The second part of the lemma says that under the same assumption, for general initial states, if we do not distinguish the local states $(1,0)$ and $(1,1)$, then the resulting process is still a Markov process. This follows from \cite[Proposition~3.1]{Swa22}.
\bl[Special choice of parameters]
Assume\label{L:special} that $X=(X^1,X^2)=(X^1_t,X^2_t)_{t\geq0}$ is a \emph{2CP(}$\lambda,\delta,\lambda_\vee,\delta_\vee,\lambda_\oplus,\delta_\oplus)$ with $\delta_\vee=\lambda_\oplus=0$. Then
\be\label{ordpres}
X^1_0(i)\geq X^2_0(i)\quad(i\in\Z^d)\quad\mbox{implies}\quad X^1_t(i)\geq X^2_t(i)\quad(i\in\Z^d,\ t\geq 0).
\ee
Moreover, if $x\mapsto\ov x$ denotes the map from $U$ to $\{(0,0),(0,1),(1,?)\}$ defined as
\be\label{forget}
\ov{(0,0)}:=(0,0),\quad\ov{(0,1)}:=(0,1),\quad\ov{(1,0)}=\ov{(1,1)}:=(1,?),
\ee
and $\ov X_t(i):=\ov{X_t(i)}$ is defined coordinatewise, then $\ov X=(\ov X_t)_{t\geq 0}$ is a Markov process.
\el


\subsection{The basic duality relation}
\label{S:basicDUAL}

Additive and cancellative duality are important tools in the study of interacting particle systems \cite{Gri79}. It is well-known that the CP is self-dual in the sense of additive systems duality. Similarly, the cCP is self-dual in the sense of cancellative systems duality. This suggests that the 2CP should also possess a self-duality. 

To present a complete picture we repeat the definitions of the additive and the cancellative duality funcion. Analogously to \cite{LS22} we define $\psi_1,\psi_2:T\times T\to T$ as
\begin{align}
\label{psi1&2}
\begin{pmatrix}
\psi_1(0,0)&\psi_1(0,1)\\
\psi_1(1,0)&\psi_1(1,1)\\
\end{pmatrix}=
\begin{pmatrix}
\psi_2(0,0)&\psi_2(0,1)\\
\psi_2(1,0)&\psi_2(1,1)\\
\end{pmatrix}=
\begin{pmatrix}
0&0\\
0&1\\
\end{pmatrix},
\end{align}
in parallel to (\ref{S_fin}) we set
\begin{align}
\label{T_fin}
\mathcal{T}_{\rm fin}:=\big\{x\in\mathcal{T}:\vert\{i\in\Z^d:x(i)\neq0\}\vert<\infty\big\}
\end{align}
and, for $x,y\in\mathcal{T}$ so that either $x\in\mathcal{T}_{\rm fin}$ or $y\in\mathcal{T}_{\rm fin}$, we define
\begin{align}\label{psib12}
\psib_1(x,y):=\bigvee_{i\in\Z^d}\psi_1\big(x(i),y(i)\big)\qquad\text{and}\qquad\psib_2(x,y):=\bigoplus_{i\in\Z^d}\psi_2\big(x(i),y(i)\big),
\end{align}
where $\vee$ and $\oplus$ are the operators defined in Section~\ref{S:PoI}. Since either $x\in\mathcal{T}_{\rm fin}$ or $y\in\mathcal{T}_{\rm fin}$, only finitely many ``summands" in the infinite ``sums" above are different from 0 and hence the expressions are well-defined. Fix $\la,\de\geq 0$, let $(X^1_t)_{t\geq 0}$ denote the CP$(\la,\de)$, and let $(X^2_t)_{t\geq 0}$ denote the cCP$(\la,\de)$.  
Then it is known \cite[Lemmas 6.6 and 6.11]{Swa22} that the contact process and cancellative contact process are self-dual in the sense that
\be\label{cancadd}
\E^x\big[\psib_i(X^i_t,y)\big]=\E^y\big[\psib_i(x,X^i_t)\big]
\qquad(x\in\mathcal{T},\ y\in\mathcal{T}_{\rm fin},\ t\geq 0,\ i=1,2),
\ee
where $\E^z$ denotes expectation with respect to the law of the process started in the initial state $X^i_0=z$ $(i=1,2,\ z\in\{x,y\})$. In general, throughout this paper, we write $\mathbb{P}^z$ and $\mathbb{E}^z$ to denote the law and expectation of a Markov process $Z=(Z_t)_{t\geq0}$ started in the initial state $Z_0=z$.

We will prove a similar self-duality for the 2CP. The first step is to find the right duality function. To this aim, we rewrite the duality functions $\psib_1,\psib_2$ in (\ref{psib12}) in such a way that the operators $\vee$ and $\oplus$ are replaced by the product in $\R$. For this purpose, we define maps $\gamma_i:T\to\R$ $(i=1,2)$ by
\begin{align}
\label{gammas}
\gamma_1(0)=1,\quad\gamma_1(1)=0\qquad\text{and}\qquad\gamma_2(0)=1,\quad\gamma_2(1)=-1.
\end{align}
Then it is easy to check that $\ga_1(x\vee y)=\ga_1(x)\cdot\ga_1(y)$ and $\ga_2(x\oplus y)=\ga_2(x)\cdot\ga_2(y)$ $(x,y\in T)$. We define, again for $x,y\in\mathcal{T}$ so that either $x\in\mathcal{T}_{\rm fin}$ or $y\in\mathcal{T}_{\rm fin}$,
\begin{align}
\label{addANDcanc}
\boldsymbol{\psi}_\text{add}(x,y):=\gamma_1\big(\psib_1(x,y)\big)\qquad\text{and}\qquad\boldsymbol{\psi}_\text{canc}(x,y):=\gamma_2\big(\psib_2(x,y)\big).
\end{align}
One then readily checks that
\begin{align*}
\boldsymbol{\psi}_\text{add}(x,y)=\prod_{i\in\Z^d}\gamma_1\big(\psi_1(x(i),y(i))\big)\qquad\text{and}\qquad\boldsymbol{\psi}_\text{canc}(x,y)=\prod_{i\in\Z^d}\gamma_2\big(\psi_2(x(i),y(i))\big),
\end{align*}
where the product is the usual product in $\R$. As $\gamma_1$ and $\gamma_2$ are bijections from $T$ to $\{0,1\}$ resp.\ to $\{-1,1\}$, (\ref{cancadd}) remains true if we replace $\psib_1$ by $\psib_{\rm add}$ and $\psib_2$ by $\psib_{\rm canc}$. Both $\psib_1$ and $\psib_\text{add}$ are known in the literature as ``the" \emph{additive duality function}. Similarly both $\psib_2$ and $\psib_\text{canc}$ are known as ``the" \emph{cancellative duality function}. With the notions introduced in Section~\ref{S:Info} below, one could more formally say that $\psib_\text{add}$ and $\psib_\text{canc}$ are (good) multiplicative representations of $\psib_1$ and $\psib_2$, respectively.

We now define, for $x=(x_1,x_2),y=(y_1,y_2)\in\mathcal{U}$ so that either $x\in\mathcal{U}_{\rm fin}$ or $y\in\mathcal{U}_{\rm fin}$,
\begin{align*}
\psib(x,y):=\psib_\text{add}(x_1,y_1)\psib_\text{canc}(x_2,y_2).
\end{align*}
One then checks that
\begin{align}
\label{psidef}
\psib(x,y)=\prod_{i\in\Z^d}\psi\big(x(i),y(i)\big),
\end{align}
where
\begin{align}
\label{frompsi1&2topsi}
\psi\big(x(i),y(i)\big)=\gamma_1\big(\psi_1(x_1(i),y_1(i))\big)\gamma_2\big(\psi_2(x_2(i),y_2(i))\big),
\end{align}
i.e.\ $\psi:U\times U\to\{-1,0,1\}$ is defined as
\begin{align}
\label{psidefloc}
\begin{pmatrix}
\psi((0,0),(0,0))&\psi((0,0),(0,1))&\psi((0,0),(1,0))&\psi((0,0),(1,1))\\
\psi((0,1),(0,0))&\psi((0,1),(0,1))&\psi((0,1),(1,0))&\psi((0,1),(1,1))\\
\psi((1,0),(0,0))&\psi((1,0),(0,1))&\psi((1,0),(1,0))&\psi((1,0),(1,1))\\
\psi((1,1),(0,0))&\psi((1,1),(0,1))&\psi((1,1),(1,0))&\psi((1,1),(1,1))\\
\end{pmatrix}
:=
\begin{pmatrix}
1&1&1&1\\
1&-1&1&-1\\
1&1&0&0\\
1&-1&0&0
\end{pmatrix}.
\end{align}
The basis of the present paper is the following duality relation.

\bp[Basic duality relation]
\label{P:dual}
For $\lambda,\delta,\lambda_\vee,\delta_\vee,\lambda_\oplus,\delta_\oplus\geq0$ let $X=(X_t)_{t\geq 0}$ and $Y=(Y_t)_{t\geq 0}$ both be a \emph{2CP(}$\lambda,\delta,\lambda_\vee,\delta_\vee,\lambda_\oplus,\delta_\oplus$\emph{)}. 
Fixing a $t\geq0$, for independent initial states such that almost surely either $X_0\in\mathcal{U}_{\rm fin}$ or $Y_0\in\mathcal{U}_{\rm fin}$, we can almost surely construct $X$ and $Y$ in such a way that for every $s\in[0,t]$ the random variables $X_s$ and $Y_{t-s}$ are independent and
\begin{align*}
[0,t]\ni s\mapsto\psib(X_s,Y_{t-s}^-)
\end{align*}
is constant, where $Y^-=(Y^-_t)_{t\geq0}$ is the c\`agl\`ad modification of $Y$, i.e.\ it is left-continuous with right limits but coincides almost everywhere with $Y$, which is c\`adl\`ag, i.e.\ right-continuous with left limits.
\ep

In fact, in the following we only need equality in expectation, i.e.\ that
\begin{align}
\label{dual}
\E\big[\psib(X_s,Y_{t-s})\big]=\E\big[\psib(X_u,Y_{t-u})\big]
\end{align}
for all $s,u\in[0,t]$. In particular, setting $s=t$ and $u=0$ and restricting ourselves to the case that $Y_0=y$ and $X_0=x$ are deterministic, this is a relation of the form (\ref{cancadd}), but with the cancellative and additive duality functions $\psib_1$ and $\psib_2$ replaced by the new duality function $\psib$. Note that, by (\ref{finpres}) and the assumption that either $X_0\in\mathcal{U}_\text{fin}$ or $Y_0\in\mathcal{U}_\text{fin}$, the expression $\psib(X_s,Y_{t-s})$ is well-defined for all $s\in[0,t]$. The following lemma highlights the strength of the duality relation (\ref{dual}).

\bl[The duality is informative]
If\label{L:info} $X$ and $X'$ are $\mathcal{U}$-valued random variables such that
\begin{align*}
\E\big[\psib(X,y)\big]=\E\big[\psib(X',y)\big]
\end{align*}
for all $y\in\mathcal{U}_{\rm fin}$, then $X$ and $X'$ are equal in distribution.
\el

We recall that a semigroup with a neutral element is called a \emph{monoid}. Examples of monoids are $(T,\vee)$, $(T,\oplus)$ and $(U,\veebar)$. These three monoids are all commutative. In the paper \cite{LS22}, the monoid $(U,\veebar)$ is called $M_{23}$, the monoid $(\{-1,0,1\},\cdot)$ is called $M_5$, and the duality function $\psi$ from (\ref{psidefloc}) is called $\psi_{235}$. In the special setting of Lemma~\ref{L:special}, the duality relation (\ref{dual}) simplifies. Assume that $\delta_\vee=\lambda_\oplus=0$ and that $X^1_0(i)\geq X^2_0(i)$ $(i\in\Z^d)$ and hence by (\ref{ordpres}) also $X^1_t(i)\geq X^2_t(i)$ $(i\in\Z^d)$ for all $t\geq 0$. Let $y\mapsto\ov y$ be the map defined in (\ref{forget}). Then (\ref{psidefloc}) and (\ref{dual}) imply that
\begin{align*}
\E\big[\ov\psib(X_s,\ov Y_{t-s})\big]=\E\big[\ov\psib(X_u,\ov Y_{t-u})\big]
\end{align*}
for all $s,u\in[0,t]$, where $\ov\psib(x,y)=\prod_{i\in\Z^d}\ov\psi\big(x(i),y(i)\big)$ with
\begin{align*}
\begin{pmatrix}
\ov\psi((0,0),(0,0))&\ov\psi((0,0),(0,1))&\ov\psi((0,0),(1,?))\\
\ov\psi((1,0),(0,0))&\ov\psi((1,0),(0,1))&\ov\psi((1,0),(1,?))\\
\ov\psi((1,1),(0,0))&\ov\psi((1,1),(0,1))&\ov\psi((1,1),(1,?))\\
\end{pmatrix}
:=
\begin{pmatrix}
1&1&1\\
1&1&0\\
1&-1&0
\end{pmatrix}.
\end{align*}
In \cite{LS22}, the duality function $\ov\psi$ is called $\psi_5$ and the sub-monoid of $U$ consisting of all $(x^1,x^2)\in U$ with $x^1\geq x^2$ is called $M_6$. The monoid $\{(0,0),(0,1),(1,?)\}$ obtained from $U$ by identifying the elements $(1,0)$ and $(1,1)$ with the single new element $(1,?)$ is isomorphic to the already mentioned $M_5$.

\subsection{Long-time behaviour}
\label{S:long-time}

We equip $\mathcal{U}$ with the product topology and define shift operators $\tet_i:\mathcal{U}\to\mathcal{U}$ by
\begin{align*}
(\tet_ix)(j):=x(j-i)\qquad\big(i,j\in\Z^d,\ x\in\mathcal{U}\big).
\end{align*}
We say that a probability measure $\mu$ on $\mathcal{T}$ or on $\mathcal{U}$ is \emph{shift-invariant} if $\mu=\mu\circ\tet_i^{-1}$ $(i\in\Z^d)$. For $x\in T$ or $x\in U$, we let $\un x$ denote the constant configuration $\un x(i):=x$ $(i\in\Z^d)$. We say that a distribution $\mu$ on $\mathcal{T}$ is \emph{non-trivial} if $\mu(\{\underline{0}\})=0$. For a distribution $\mu^\prime$ on $\mathcal{U}$ we define non-triviality to hold if
\begin{align*}
\mu^\prime\big(\big\{\underline{(0,0)}\big\}\big)=0.
\end{align*}

It is well-known \cite[Theorem~6.35]{Swa22} that the CP($\lambda,\delta$) with $\lambda+\delta>0$ started in a non-trivial shift-invariant distribution converges weakly to a (time-) invariant distribution $\bar{\nu}$ called the \emph{upper invariant law} of the contact process. Similarly, it is known \cite[Theorem~1.2\ \&\ Theorem~1.3]{BDD91} that the cCP($\lambda,\delta$) with $\lambda+\delta>0$ started in a non-trivial shift-invariant distribution converges weakly to an invariant distribution $\dot{\nu}$, that we call, in accordance with \cite{SS08a}, the \emph{odd upper invariant law} of the cancellative contact process. 

Letting $\delta_{\underline{0}}$ denote the Dirac measure concentrated on the ``all 0" configuration $\underline{0}$, $\bar{\nu}$ and $\dot{\nu}$ may or may not differ from $\delta_{\underline{0}}$ depending on the choice of the model parameters $\lambda$ and $\delta$. For a CP($\lambda,\delta$) ($\lambda+\delta>0$) there exists a critical value $\lambda_{\text{CP}}=\lambda_{\text{CP}}(d)\in(0,\infty)$ (dependent on the dimension $d$) such that $\bar{\nu}\neq\delta_{\underline{0}}$ if and only if $\lambda/\delta>\lambda_\text{CP}$ \cite{BG90}. Here and in the following we set $x/0=\infty$ for $x\in(0,\infty)$. For the cCP we can define $\lambda_{\text{cCP}}^\pm=\lambda_{\text{cCP}}^\pm(d)$ as
\begin{align*}
\lambda_\text{cCP}^-&:=\inf\{\lambda\geq0:\text{the odd upper invariant law of the cCP}(\lambda,1)\ \text{does not equal}\ \delta_{\underline{0}}\},\\
\lambda_\text{cCP}^+&:=\sup\{\lambda\geq0:\text{the odd upper invariant law of the cCP}(\lambda,1)\ \text{equals}\ \delta_{\underline{0}}\}.
\end{align*}
It is known that $\lambda_\text{cCP}^+<\infty$ (\cite[Theorem~1.1]{BDD91} \& Proposition~\ref{P:surv} below) and considering the 2CP started with no $i\in\mathbb{Z}^d$ in the local state $(0,1)$, i.e.\ the coupling from Section~\ref{S:PoI}, shows that $\lambda_\text{CP}\leq\lambda^-_\text{cCP}$, thus
\begin{align*}
0<\lambda_\text{CP}\leq\lambda^-_\text{cCP}\leq\lambda^+_\text{cCP}<\infty.
\end{align*}
Simulations suggest that $\lambda^-_\text{cCP}=\lambda^+_\text{cCP}$ in all dimensions but proving this is a long-standing open problem that due to the non-monotone nature of the process seems very difficult. Moreover, one would assume that actually $\lambda_\text{CP}<\lambda^-_\text{cCP}$ holds. In dimension one, using the bound $\lambda_\text{CP}(1)\leq1.942$, proved in \cite{Lig95}, and the following result we can prove just that. 
\begin{proposition}[Lower bound for $\mathbf{\lambda_\text{cCP}^-(1)}$]
\label{P:lambda_cCP(1)<=2}
One has $\lambda^-_\text{\emph{cCP}}(1)\geq2$.
\end{proposition}
The proof of Proposition~\ref{P:lambda_cCP(1)<=2}, which can be found in Section~\ref{S:surv}, goes in two steps. First, cancellative duality is used to show that one may equivalently prove that the cCP$(\la,\de)$ with $\la<2\de$ started from a finite initial state dies out almost surely. Next, Dynkin's formula is used to prove that for $\la\leq 2\de$, the position of the right-most (left-most) one is a supermartingale (submartingale), which is shown to imply extinction of finite processes. This last argument is adapted from \cite{Sud98} who uses it to prove for the CP the bound $\lambda_\text{CP}\geq1$ and shows that better bounds can be obtained by a refined argument that constructs a supermartingale that depends not only on the position of the right-most one but also takes into consideration the configuration on finitely many sites left of the right-most one. As these methods are essentially one-dimensional in nature, it is not clear how to generalise Proposition~\ref{P:lambda_cCP(1)<=2} to higher dimensions.


In this paper we are interested in the long-time behaviour of the 2CP started in a shift-invariant distribution. We set
\begin{align*}
\mathcal{U}_{(0,\ast)}&:=\big\{x=(x_1,x_2)\in\mathcal{U}:x_1=\un0\big\},\\
\mathcal{U}_{(\ast,0)}&:=\big\{x=(x_1,x_2)\in\mathcal{U}:x_2=\un0\big\},\\
\mathcal{U}_\text{mix}&:=\mathcal{U}\setminus(\mathcal{U}_{(0,\ast)}\cup\mathcal{U}_{(\ast,0)}).
\end{align*}
The known results for CP and cCP imply that the 2CP $X=(X_t)_{t\geq0}=(X^1_t,X^2_t)_{t\geq0}$ started in a non-trivial shift-invariant distribution on $\mathcal{U}_{(\ast,0)}$ converges weakly to $\bar{\nu}\otimes\delta_{\underline{0}}$. Analogously, the 2CP started in a non-trivial shift-invariant distribution on $\mathcal{U}_{(0,\ast)}$ converges weakly to $\delta_{\underline{0}}\otimes\dot{\nu}$. If $X$ is started in a non-trivial shift-invariant distribution on $\mathcal{U}_\text{mix}$, then the laws of $X^1_t$ and $X^2_t$ individually converge weakly as $t\to\infty$ to $\bar{\nu}$ and $\dot{\nu}$, respectively. However, as a measure on a product space is in general not determined by its marginals, the long-time behaviour of the joint law of $X_t=(X^1_t,X^2_t)$ is less straightforward.
A priori there might, for example, exist an increasing sequence $(t_n)_{n\in\mathbb{N}}$ so that the sequence of laws of $(X_{t_n})_{n\in\mathbb{N}}$ has several cluster points all having the marginal distributions $\bar{\nu}$ and $\dot{\nu}$, respectively. Or the law of $X$ might converge weakly to different distributions depending on where on $\mathcal{U}_\text{mix}$ its initial law is supported. We will use the duality function $\boldsymbol{\psi}$ to show that the behaviour outlined in the last two sentences does not occur.
\begin{theorem}[Joint invariant law]
\label{T:invlaw}
Let $X=(X^1,X^2)=(X^1_t,X^2_t)_{t\geq0}$ be a \emph{2CP} with parameters $\lambda,\delta,\lambda_\vee,\delta_\vee,\lambda_\oplus,\delta_\oplus\geq0$ so that $\lambda+\lambda_\vee+\delta+\delta_\vee>0$ and $\lambda+\lambda_\oplus+\delta+\delta_\oplus>0$. Then $X$ has an invariant law $\nu$ that is uniquely characterised by the relation
\begin{align*}
\int\psib(x,y)\ \di\nu(x)=\P^y\left[\exists t\geq 0:X_t=\underline{(0,0)}\right]
\qquad\big(y\in\mathcal{U}_{\rm fin}\big).
\end{align*}
If $X$ is started in a shift-invariant initial law that is concentrated on $\mathcal{U}_\text{\emph{mix}}$, then
\begin{align}
\label{nuconv}
\P\big[X_t\in\,\cdot\,\big]\Asto{t}\nu.
\end{align}
\end{theorem}
As usual, the symbol $\Rightarrow$ in (\ref{nuconv}) denotes weak convergence. Note that (\ref{nuconv}) implies that $\nu$ is (as $\bar{\nu}$ and $\dot{\nu}$) shift-invariant. In the special case that $\delta_\vee=\lambda_\oplus=0$, corresponding to the monotone coupling of CP and cCP, one has that
\begin{align*}
\nu\big(\big\{x\in\mathcal{U}:\exists i\in\mathbb{Z}^d:x(i)=(0,1)\big\}\big)=0,
\end{align*}
as we can chose a shift-invariant initial law that is concentrated on $\mathcal{U}_\text{mix}$ with the above property. This property is then preserved by the dynamics. One example of such an initial law would be the Dirac measure concentrated on $\underline{(1,1)}$. Thus, as long as the initial distribution of this special 2CP is shift-invariant and concentrated on $\mathcal{U}_\text{mix}$, the law of this 2CP converges weakly to a monotonically coupled law, no matter how high the density of $(0,1)$s was in the initial distribution. 

Taking into account our earlier remarks about initial laws on $\mathcal{U}_{(0,\ast)}$ and  $\mathcal{U}_{(\ast,0)}$, one can conclude (compare \cite[Corollary~6.39]{Swa22}) that all shift-invariant invariant laws of the 2CP are convex combinations of $\delta_{\underline{0}}\otimes\delta_{\underline{0}}$, $\bar{\nu}\otimes\delta_{\underline{0}}$, $\delta_{\underline{0}}\otimes\dot{\nu}$ and $\nu$.

\subsection{Outline}

The paper is structured as follows. In Section \ref{S:MProc} we provide a proof for Proposition~\ref{P:dual}. If fact, we prove in Theorem~\ref{T:dual} a generalisation of Proposition~\ref{P:dual} that is independent of our process of interest, so that it can directly be applied to further processes. Section~\ref{S:Info} deals with the proof of Lemma~\ref{L:info}. Also here we prove in Proposition~\ref{P:genINFO} a generalisation of Lemma~\ref{L:info}. Additionally, towards the goal of classifying the dualities found in \cite{LS22} regarding their ability to determine laws of processes uniquely, we introduce two notions and show that they basically coincide in our setup. In Section~\ref{S:weak} we prove Theorem~\ref{T:invlaw}. As Proposition~\ref{P:lambda_cCP(1)<=2} is independent of the monoid dualities from \cite{LS22}, we prove it last. Its proof is found in Section~\ref{S:surv}. Finally, in Appendix~\ref{S:App} we show how Lemma~\ref{L:indeter}, an auxiliary result we use for the proof of Theorem~\ref{T:dual}, follows from a corollary from \cite{SS08a}. As this corollary is stated in \cite{SS08a} in a rather general form, we decided to repeat the definitions from \cite{SS08a}, slightly reformulate the result and move this discussion to the appendix.


\section{Monoid duality for interacting particle systems}
\label{S:MProc}

In \cite{LS22} a duality theory is developed for Markov processes with state space of the form $S^\La$ where $S$ is a finite commutative monoid and $\La$ is a finite set. Here we generalise this to countable $\La$ which allows us to define duality relations for interacting particle systems. For the special cases of additive and cancellative dualities infinite $\La$ have already been treated in \cite[Chapter~6.6 \& Chapter~6.7]{Swa22}.

We start by extending the concept of duality between monoids (i.e.\ semigroups with a neutral element) presented in \cite{LS22} to monoids that carry a topology. We say that a monoid $(M,+)$ is a \emph{topological monoid} if it is equipped with a topology so that the map $M\times M\ni (x,y)\mapsto x+y\in M$ is continuous, where $M\times M$ is equipped with the product topology. For a second topological monoid $(N,+)$ we denote by $\mathcal{H}(M,N)$ the space of all continuous monoid homomorphisms, i.e.\ continuous functions from $M$ to $N$ that preserve the operation and map the neutral element of $M$ to the neutral element of $N$. Throughout this paper we always equip finite and countable monoids with the discrete topology, so that every finite or countable monoid is a topological monoid. This makes every function between two finite or countable monoids continuous. Thus, if $N$ and $M$ are finite, the space $\mathcal{H}(M,N)$ defined above coincides with the space of all monoid homomorphisms (called $\mathcal{H}(M,N)$ in \cite{LS22}).

Let $M_1$, $M_2$ and $N$ be topological monoids. We say that $M_1$ is \emph{$N$-dual} to $M_2$ with \emph{duality function} $\psi$ if the following conditions are satisfied:
\begin{enumerate}
\item $\psi(x_1,y)=\psi(x_2,y)$ for all $y\in M_2$ implies $x_1=x_2$ $(x_1,x_2\in M_1)$,
\item $\mathcal{H}(M_1,N)=\{\psi(\,\cdot\,,y):y\in N\}$,
\item $\psi(x,y_1)=\psi(x,y_2)$ for all $x\in M_1$ implies $y_1=y_2$ $(y_1,y_2\in M_2)$,
\item $\mathcal{H}(M_2,N)=\{\psi(x,\,\cdot\,):x\in M_1\}$.
\end{enumerate}
As we equip finite monoids with the discrete topology, the definition above coincides with the definition of duality between monoids from \cite{LS22} if $M_1,M_2$ and $N$ are finite.

Repeating the definition from \cite{LS22}, for arbitrary spaces $\mathcal{X},\mathcal{Y}$ and $\mathcal{Z}$ we say that the map $m:\mathcal{X}\to\mathcal{X}$ is \emph{dual} to the map $\hat{m}:\mathcal{Y}\to\mathcal{Y}$ with respect to the \emph{duality function} $\psi:\mathcal{X}\times\mathcal{Y}\to\mathcal{Z}$ if
\begin{align*}
\psi(m(x),y)=\psi(x,\hat{m}(y))\qquad(x\in\mathcal{X},\ y\in\mathcal{Y}).
\end{align*}

In parallel to \cite{Swa22} we say that a map $m:\mathcal{X}\to\mathcal{X}$ \emph{preserves} a set $\mathcal{H}$ of functions from $\mathcal{X}$ to $\mathcal{Y}$ if
\begin{align*}
f\circ m\in\mathcal{H}\quad\text{for all}\quad f\in\mathcal{H}.
\end{align*}
The following proposition is the analogue of \cite[Proposition~5]{LS22}, that is formulated for dualities between monoids without attached topologies.

\begin{proposition}[Maps having a dual]
\label{P:mapsWduals}
Let $S$, $R$ and $T$ be commutative topological monoids such that $S$ is $T$-dual to $R$ with duality function $\psi$. Then a map $m:S\to S$ has a dual map $\hat{m}:R\to R$ with respect to $\psi$ if and only if $m$ preserves $\mathcal{H}(S,T)$. The dual map $\hat{m}$, if it exists, is unique and preserves $\mathcal{H}(R,T)$.
\end{proposition}
\begin{proof}
If $m:S\to S$ preserves $\mathcal{H}(S,T)$, then, by property (ii) of the definition of duality, for all $y\in R$ one has $\psi(m(\,\cdot\,),y)\in\mathcal{H}(S,T)$. Applying property (ii) again, it follows that there exists an $\hat{m}(y)\in R$ such that $\psi(m(\,\cdot\,),y)=\psi(\,\cdot\,,\hat{m}(y))$. Property (iii) of the definition of duality implies that $\hat{m}(y)$ is unique. This shows that $m$ has a unique dual map $\hat{m}:R\to R$ if $m$ preserves $\mathcal{H}(S,T)$.

On the other hand, if $m:S\to S$ has a dual map $\hat{m}:R\to R$, then $\psi(m(\,\cdot\,),y)=\psi(\,\cdot\,,\hat{m}(y))$, i.e.\ $m$ preserves $\{\psi(\,\cdot\,,y):y\in R\}$. By property (ii) of the definition of duality $m$ then also preserves $\mathcal{H}(S,T)$. This finishes the proof that $m:S\to S$ has a dual map $\hat{m}:R\to R$ if and only if $m$ preserves $\mathcal{H}(S,T)$.

Finally, if $\hat{m}$ exists, then it has $m$ as a dual map with respect to $\psi^\dagger:R\times S\to T$ defined as $\psi^\dagger(y,x):=\psi(x,y)$  $(y\in R,\ x\in S)$, and the previously proved statement implies that $\hat{m}$ has to preserve $\mathcal{H}(R,T)$.
\end{proof}

Clearly, any $m\in\mathcal{H}(S,S)$ preserves $\mathcal{H}(S,T)$. Conversely, if the assumptions on $S$, $T$ and $R$ from Proposition~\ref{P:mapsWduals} are satisfied and $m:S\to S$ preserves $\mathcal{H}(S,T)$, then the proof of \cite[Proposition~5]{LS22} shows that $m:S\to S$ has to be a monoid homomorphism. However, while duality implies that $\psi(m(\,\cdot\,),y)$ is continuous for all $y\in R$, we do not know if $m$ itself always has to be continuous. 



We are especially interested in countable products of topological monoids as we will view state spaces of an interacting particle system as such products. Let, throughout this section, $\Lambda$ be a countable set. For a topological monoid $M$ with $\vert M\vert\geq 2$ we equip $M^\Lambda$ with the product topology, making this uncountable monoid a topological monoid. We define the countable sub-monoid $M^\Lambda_\text{fin}\subset M^\Lambda$ as
\begin{align*}
M^\Lambda_\text{fin}:=\{x\in M^\Lambda:\vert\{i\in\Lambda:x(i)\neq0\}\vert<\infty\},
\end{align*}
where 0 denotes the neutral element of $M$. As in Section~\ref{S:long-time}, we denote by $\un 0$ the constant configuration with $\un 0(i)=0$ for all $i\in\Lambda$ that is the neutral element of $M^\Lambda_\text{fin}$ and $M^\Lambda$.

Before we investigate duality between such ``product monoids" we collect some definitions and results for general product spaces from \cite{Swa22} that we will need in the following. Let $L$ and $V$ be arbitrary spaces. For a function $f:L^\Lambda\to V$ we say that $j\in\Lambda$ is \emph{f-relevant} if
\begin{align*}
\exists x_1,x_2\in L^\Lambda:f(x_1)\neq f(x_2)\ \text{but}\ x_1(k)=x_2(k)\ \forall k\neq j.
\end{align*}
We set
\begin{align*}
\mathcal{R}(f):=\{j\in\Lambda:j\ \text{is}\ f\text{-relevant}\}
\end{align*}
and cite the following result \cite[Lemma~4.13]{Swa22}. 
\begin{lemma}[Continuous maps]
\label{L:CONTmapPRODtop}
Let $L$ and $V$ be finite sets equipped with the discrete topology. A map $f:L^\Lambda\to V$ is continuous with respect to the product topology if and only if the following two conditions hold:
\begin{itemize}
\item[\emph{(i)}]$\mathcal{R}(f)$ is finite.
\item[\emph{(ii)}]If $x_1,x_2\in L^\Lambda$ satisfy $x_1(j)=x_2(j)$ for all $j\in\mathcal{R}(f)$, then $f(x_1)=f(x_2)$.
\end{itemize}
\end{lemma}
Let $L$ be finite. For any map $\mathtt{m}:L^\Lambda\to L^\Lambda$ and $i\in\Lambda$ we define $\mathtt{m}[i]:L^\Lambda\to L$ as
\begin{align*}
\mathtt{m}[i](x):=\mathtt{m}(x)(i)\qquad(x\in L^\Lambda).
\end{align*}
Moreover, we let
\begin{align*}
\mathcal{D}(\mathtt{m}):=\{i\in\Lambda:\exists x\in L^\Lambda: \mathtt{m}[i](x)\neq x(i)\}.
\end{align*}
We say that a map $\mathtt{m}:L^\Lambda\to L^\Lambda$ is \emph{local} if
\begin{align*}
\text{(i)}\ \mathtt{m}\ \text{is continuous}\qquad\text{and}\qquad\text{(ii)}\ \mathcal{D}(\mathtt{m})\ \text{is finite}.
\end{align*}
For a finite monoid $M$ we denote by $\mathcal{H}_\text{loc}(M^\Lambda,M^\Lambda)$ the space of all maps $\mathtt{m}\in\mathcal{H}(M^\Lambda,M^\Lambda)$ that are local. As we equip, according to our conventions, $M^\Lambda$ with the product topology and $M$ with the discrete one, $\mathtt{m}:M^\Lambda\to M^\Lambda$ is local if and only if $\mathcal{D}(\mathtt{m})$ is finite and $\mathtt{m}[j]$ satisfies the conditions of Lemma~\ref{L:CONTmapPRODtop} for all $j\in\Lambda$. Note that every $\mathtt{m}\in\mathcal{H}_\text{loc}(M^\Lambda,M^\Lambda)$ maps $M^\Lambda_\text{fin}$ into itself. 

Let throughout the rest of this section $(S,\odot)$, $(R,\raisebox{-0.35mm}{$\boxdot$})$ and $(T,\otimes)$ be commutative finite (and hence topological) monoids and assume that $S$ is $T$-dual to $R$ with duality function $\psi:S\times R\to T$. We denote all three neutral elements by 0 and define $\Psi:S^\Lambda\times R^\Lambda_\text{fin}\to T$ by
\begin{align}
\label{Psi}
\Psi(\mathbf{x},\mathbf{y}):=\bigotimes_{i\in\Lambda}\psi\big(\mathbf{x}(i),\mathbf{y}(i)\big)\qquad\big(\mathbf{x}\in S^\Lambda,\ \mathbf{y}\in R^\Lambda_\text{fin}\big).
\end{align}
Note that $\Psi$ is well-defined as for all but finitely many $i\in\Lambda$ one has $y(i)=0$ and $\psi(\,\cdot\,,0)=o$ due to property (iv) of the definition of duality, where $o:S\to T$ is the function that is constantly 0. In general, for all monoids $M$ and $N$, let $\text{id}\in\mathcal{H}(M,M)$ denote the identity and $o\in\mathcal{H}(M,N)$ the function constantly 0. Using Lemma~\ref{L:CONTmapPRODtop} we can prove the following.
\begin{proposition}[Duality on product spaces]
\label{P:topdualSR}
Let $S,R,T$ be finite commutative monoids. If $S$ is $T$-dual to $R$ with duality function $\psi$, then $S^\Lambda$ is $T$-dual to $R^\Lambda_\text{\emph{fin}}$ with duality function $\Psi$.
\end{proposition}
\begin{proof}
The properties (i) and (iii) from the definition of duality follow directly from the corresponding properties of the duality between $S$ and $R$. To be more precise, assuming that $\Psi(\mathbf{x}_1,\mathbf{y})=\Psi(\mathbf{x}_2,\mathbf{y})$ for all $\mathbf{y}\in R^\Lambda_\text{fin}$ in particular implies for $i\in\Lambda$ and $y\in R$ that
\begin{align*}
\psi(\mathbf{x}_1(i),y)=\Psi(\mathbf{x}_1,y^i)=\Psi(\mathbf{x}_2,y^i)=\psi(\mathbf{x}_2(i),y),
\end{align*}
where $y^i\in R^\Lambda_\text{fin}$ is defined as
\begin{align}
\label{y^i}
y^i(j)=\begin{cases}
y&\text{if}\ j=i,\\
0&\text{else},\\
\end{cases}\qquad(j\in\Lambda).
\end{align}
Hence, the fact that $S$ is $T$-dual to $R$ implies that $\mathbf{x}_1(i)=\mathbf{x}_2(i)$ for all $i\in\Lambda$ and thus $\mathbf{x}_1=\mathbf{x}_2$. Property (iii) follows in the same way.

The fact that $\Psi(\,\cdot\,,\mathbf{y})$ and $\Psi(\mathbf{x},\,\cdot\,)$ are monoid homomophisms for all $\mathbf{y}\in R^\Lambda_\text{fin}$ and for all $\mathbf{x}\in S^\Lambda$, respectively, also follows directly from the duality between $S$ and $R$ and the definition of $\Psi$. As $R^\Lambda_\text{fin}$ is countable this implies $\Psi(\mathbf{x},\,\cdot\,)\in\mathcal{H}(R^\Lambda_\text{fin},T)$. For $\mathbf{y}\in R^\Lambda_\text{fin}$ we have that $\mathcal{R}\big(\Psi(\,\cdot\,,\mathbf{y})\big)=\{j\in\Lambda:\mathbf{y}(j)\neq0\}$, so $\Psi(\,\cdot\,,\mathbf{y})$ satisfies the conditions of Lemma~\ref{L:CONTmapPRODtop} and hence also $\Psi(\,\cdot\,,\mathbf{y})\in\mathcal{H}(S^\Lambda,T)$.

To prove the implication $\subset$ in property~(iv) from the definition of duality, assume that $g\in\mathcal{H}(R^\Lambda_\text{fin},T)$. Then using (\ref{y^i}), for each $i\in\Lambda$, we define $g_i:R\to T$ as $g_i(y):=g(y^i)$ $(i\in\Lambda)$. The fact that $g\in\mathcal{H}(R^\Lambda_\text{fin},T)$ directly implies that $g_i\in\mathcal{H}(R,T)$, and the duality between $S$ and $R$ implies that there exists an $x_i\in S$ such that $g_i=\psi\big(x_i,\,\cdot\,\big)$. Defining $\mathbf{x}\in S^\Lambda$ by $\mathbf{x}(i):=x_i$, one has for $\mathbf{y}\in R^\Lambda_\text{fin}$ that
\begin{align*}
g(\mathbf{y})&=g\left(\bigboxdot\nolimits_{i:\mathbf{y}(i)\neq0}\mathbf{y}(i)^i\right)=\bigotimes_{i:\mathbf{y}(i)\neq0}g\big(\mathbf{y}(i)^i\big)=\bigotimes_{i:\mathbf{y}(i)\neq0}g_i\big(\mathbf{y}(i)\big)=\bigotimes_{i:\mathbf{y}(i)\neq0}\psi\big(x_i,\mathbf{y}(i)\big)\\&=\Psi(\mathbf{x},\mathbf{y}),
\end{align*}
which finishes the proof of property~(iv) from the definition of duality.

Lastly, we prove the implication $\subset$ in property~(ii) from the definition of duality. We assume that $f\in\mathcal{H}(S^\Lambda,T)$. Then Lemma~\ref{L:CONTmapPRODtop} implies that there exists a finite set $\Delta\subset\Lambda$ such that $f$ only depends on the coordinates in $\Delta$. Letting for $\mathbf{x}\in S^\Lambda$ the restriction $\mathbf{x}_\Gamma$ to some set $\Gamma\subset\Lambda$ be defined as
\begin{align*}
\mathbf{x}_\Gamma(j):=\begin{cases}
\mathbf{x}(j)&\text{if}\ j\in\Gamma,\\
0&\text{else},\\
\end{cases}\qquad(j\in\Lambda),
\end{align*}
we see that
\begin{align*}
f(\mathbf{x})=f\big(\mathbf{x}_{\Delta^\mathrm{c}}\odot\mathbf{x}_{\Delta}\big)=f\big(\mathbf{x}_{\Delta^\mathrm{c}}\big)\otimes\bigotimes_{i\in\Delta}f\big(\mathbf{x}(i)^i\big),
\end{align*}
where $x^i\in S^\Lambda_\text{fin}$ is defined as $y^i\in R^\Lambda_\text{fin}$ in (\ref{y^i}). But as $f$ does not depend on $\Delta^\mathrm{c}$ we conclude that
\begin{align*}
f\big(\mathbf{x}_{\Delta^\mathrm{c}}\big)=f\big(\underline{0}_{\Delta^\mathrm{c}}\big)=f(\underline{0})=0.
\end{align*}
Analogously to above we can now define $\mathbf{y}\in R^\La_\text{fin}$ by $\mathbf{y}(i):=y_i$ for $i\in\De$ and $\mathbf{y}(i):=0$ for $i\in\De^\mathrm{c}$, where $y_i\in R$ satisfies $f\big(\mathbf{x}(i)^i\big)=\psi(\mathbf{x}(i),y_i)$  independent of the value of $\mathbf{x}(i)$. 
Then $f=\Psi\big(\,\cdot\,,\mathbf{y}\big)$, which finishes the proof of property~(ii) from the definition of duality and thus the proof is complete.
\end{proof}

Having proved the duality between $S^\Lambda$ and $R^\Lambda_\text{fin}$, Proposition~\ref{P:mapsWduals} and the remarks below it imply that every $m\in\mathcal{H}(S^\Lambda,S^\Lambda)$ has a unique dual map with respect to $\Psi$. In fact, using the definition of duality and the properties of the product topology it is easy to see that $m:S^\Lambda\to S^\Lambda$ has a unique dual map with respect to $\Psi$ if and only if $m\in\mathcal{H}(S^\Lambda,S^\Lambda)$.

However, it is not clear how to compute the dual map of $m\in\mathcal{H}(S^\Lambda,S^\Lambda)$ in general, so we will focus on local monoid homomorphisms, for which we will be able to compute the dual maps explicitly. The following lemma generalises \cite[Lemma~7]{LS22} to infinite $\Lambda$. 
\begin{lemma}[Local monoid homomorphisms]
\label{L:L_HwithMAT}
Let $(S,\odot)$ be a finite monoid. Let $M=(M_{ij})_{i,j\in\Lambda}$ be an infinite matrix with values in $\mathcal{H}(S,S)$ such that the set
\begin{align}
\label{EnotID}
\Delta:=\big\{(i,j)\in\Lambda^2:i\neq j,\ M_{ij}\neq o\big\}\cup\big\{(i,i)\in\Lambda^2: M_{ii}\neq\text{\emph{id}}\big\}
\end{align}
is finite. Then setting
\begin{align}
\label{DefMATL_H}
\mathtt{m}[j](\mathbf{x}):=\bigodot_{i\in\Lambda}M_{ij}\big(\mathbf{x}(i)\big)\qquad\big(j\in\Lambda,\ \mathbf{x}\in S^\Lambda\big)
\end{align}
defines a map $\mathtt{m}\in\mathcal{H}_\text{\emph{loc}}(S^\Lambda,S^\Lambda)$. Conversely, each $\mathtt{m}\in\mathcal{H}_\text{\emph{loc}}(S^\Lambda,S^\Lambda)$ is of this form.
\end{lemma}
\begin{proof}
First assume that $\mathtt{m}$ is of the form (\ref{DefMATL_H}). Then $\mathtt{m}$ is well-defined as $\Delta$ from (\ref{EnotID}) is finite. As $M$ takes values in $\mathcal{H}(S,S)$ it follows readily that $\mathtt{m}[j]\in\mathcal{H}(S^\Lambda,S)$ for all $j\in\Lambda$, thus $\mathtt{m}\in\mathcal{H}(S^\Lambda,S^\Lambda)$. Let $j\in\Lambda$. One sees that
\begin{align*}
\mathcal{R}(\mathtt{m}[j])=\begin{cases}
\{i\in\Lambda\setminus\{j\}:(i,j)\in\Delta\}\cup\{j\}&\text{if}\ M_{jj}\neq o,\\
\{i\in\Lambda\setminus\{j\}:(i,j)\in\Delta\}&\text{if}\ M_{jj}=o.\\
\end{cases}
\end{align*}
In both cases $\mathcal{R}(\mathtt{m}[j])$ satisfies the conditions of Lemma~\ref{L:CONTmapPRODtop}. Additionally
\begin{align*}
\mathcal{D}(\mathtt{m})=\{j\in\Lambda:\exists i\in\Lambda:(i,j)\in\Delta\}
\end{align*}
is finite and it follows that $\mathtt{m}$ is local, so $\mathtt{m}\in\mathcal{H}_\text{loc}(S^\Lambda,S^\Lambda)$.

Now assume that $\mathtt{m}\in\mathcal{H}_\text{loc}(S^\Lambda,S^\Lambda)$. In particular, one has that $\mathtt{m}[j]:S^\Lambda\to S$ is continuous for all $j\in\Lambda$ by the properties of the product topology. Moreover, $\mathcal{D}(\mathtt{m})\subset\Lambda$ is finite and, by definition, for $j\in\mathcal{D}(\mathtt{m})^\mathrm{c}$ one has $\mathtt{m}[j](\mathbf{x})=\mathbf{x}(j)$ for all $\mathbf{x}\in S^\Lambda$. Due to Lemma~\ref{L:CONTmapPRODtop}, for each $j\in\mathcal{D}(\mathtt{m})$ the set $\mathcal{R}(\mathtt{m}[j])$ is finite and we can identify $\mathtt{m}[j]$ with a map $\mathtt{m}[j]\vert_{\mathcal{R}(\mathtt{m}[j])}:S^{\mathcal{R}(\mathtt{m}[j])}\to S$. By \cite[Lemma~7]{LS22} there exists a vector $M^j=(M^j_i)_{i\in\mathcal{R}(\mathtt{m}[j])}$ with coordinates in $\mathcal{H}(S,S)$ such that
\begin{align*}
\mathtt{m}[j]\vert_{\mathcal{R}(\mathtt{m}[j])}(\mathbf{x})=\bigodot_{i\in\mathcal{R}(\mathtt{m}[j])}M^j_i\big(\mathbf{x}(i)\big)\qquad\big(\mathbf{x}\in S^{\mathcal{R}(\mathtt{m}[j])}\big).
\end{align*}
Defining now $M=(M_{ij})_{i,j\in\Lambda}$ as
\begin{align*}
M_{ij}:=\begin{cases}
M_i^j&\text{if}\ j\in\mathcal{D}(\mathtt{m}),i\in\mathcal{R}(\mathtt{m}[j]),\\
\text{id}&\text{if}\ i=j\notin\mathcal{D}(\mathtt{m}),\\
o&\text{else},
\end{cases}
\end{align*}
gives a representation of $\mathtt{m}[j]$ for all $j\in\Lambda$ as in (\ref{DefMATL_H}) with the property that the set $\Delta$ from (\ref{EnotID}) is finite. This completes the proof.
\end{proof}

As already claimed, with the help of the above lemma we can compute the dual function of each $m\in\mathcal{H}_\text{loc}(S^\Lambda,S^\Lambda)$.
\begin{proposition}[Dual local homomorphisms]
\label{P:dualMapL}
\hspace{-3.5pt} Let $S, R, T$ be finite commutative monoids so that $S$ is $T$-dual to $R$ with duality function $\psi$. For each $\mathtt{m}\in\mathcal{H}_\text{\emph{loc}}(S^\Lambda,S^\Lambda)$ there exists a map $\hat{\mathtt{m}}\in\mathcal{H}_\text{\emph{loc}}(R^\Lambda,R^\Lambda)$ so that the restriction of $\hat{\mathtt{m}}$ to $R^\Lambda_\text{\emph{fin}}$ is the unique dual map of $\mathtt{m}$ with respect to the duality function $\Psi$ from (\ref{Psi}). If $M=(M_{ij})_{i,j\in\Lambda}$ denotes the matrix from Lemma~\ref{L:L_HwithMAT} such that (\ref{DefMATL_H}) holds, then $\hat{\mathtt{m}}$ is given via
\begin{align}
\label{hatM}
\hat{\mathtt{m}}[j](\mathbf{y})=\bigboxdot_{j\in\Lambda}\widehat{M}_{ij}\big(\mathbf{y}(j)\big)\qquad(j\in\Lambda,\,\mathbf{y}\in R^\Lambda),
\end{align}
where, for $i,j\in\Lambda$, $\widehat{M}_{ij}\in\mathcal{H}(R,R)$ is the (unique) dual map of $M_{ij}\in\mathcal{H}(S,S)$ with respect to the duality function $\psi$.
\end{proposition}
\begin{proof}
Let $\mathbf{x}\in S^\Lambda$, $\mathbf{y}\in R^\Lambda_\text{fin}$ and let $\hat{\mathtt{m}}$ be defined via (\ref{hatM}). Note that $\hat{\mathtt{m}}$ indeed maps $R^\Lambda_\text{fin}$ into itself as $\Delta$ from (\ref{EnotID}) is finite for $\mathtt{m}$ and the (unique) dual maps of $o,\text{id}\in\mathcal{H}(S,S)$ with respect to $\psi$ are $o\in\mathcal{H}(R,R)$ and $\text{id}\in\mathcal{H}(R,R)$, respectively. Moreover, Lemma~\ref{L:L_HwithMAT} implies that $\hat{\mathtt{m}}\in\mathcal{H}_\text{loc}(R^\Lambda,R^\Lambda)$. We compute that
\begin{align*}
\Psi(\mathtt{m}(\mathbf{x}),\mathbf{y})&=\bigotimes_{j\in\Lambda}\psi\left(\bigodot\nolimits_{i\in\Lambda}M_{ij}(\mathbf{x}(i)),\,\mathbf{y}(j)\right)=\bigotimes_{i,j\in\Lambda}\psi\big(M_{ij}(\mathbf{x}(i)),\,\mathbf{y}(j)\big)\\
&=\bigotimes_{i,j\in\Lambda}\psi\big(\mathbf{x}(i),\,\widehat{M}_{ij}(\mathbf{y}(j))\big)=\bigotimes_{i\in\Lambda}\psi\left(\mathbf{x}(i),\bigboxdot\nolimits_{j\in\Lambda}\widehat{M}_{ij}(\mathbf{y}(j))\right)\\
&=\Psi(\mathbf{x},\hat{\mathtt{m}}(\mathbf{y})).
\end{align*}
Uniqueness of the dual map follows directly from property (iii) of the duality between $S^\Lambda$ and $R^\Lambda_\text{fin}$ established in Proposition~\ref{P:topdualSR}.
\end{proof}

%
%

We are now ready to apply the non-probabilistic results above to Markov processes. Let $S,R$ and $T$ still be the finite monoids from above and let $\mathcal{G}$ be a countable collection of maps in $\mathcal{H}_\text{loc}(S^\Lambda,S^\Lambda)$. We are considering two formal Markov generators $G$ and $\widehat{G}$ defined as
\begin{align}
\label{monG}
Gf(\mathbf{x})&:=\sum_{\mathtt{m}\in\mathcal{G}}r_\mathtt{m}\big(f(\mathtt{m}(\mathbf{x}))-f(\mathbf{x})\big)\qquad(\mathbf{x}\in S^\Lambda),
\intertext{and}
\label{monGhat}
\widehat{G}g(\mathbf{y})&:=\sum_{\mathtt{m}\in\mathcal{G}}r_\mathtt{m}\big(g(\hat{\mathtt{m}}(\mathbf{y}))-g(\mathbf{y})\big)\qquad(\mathbf{y}\in R^\Lambda_\text{fin}),
\end{align}
where $\hat{\mathtt{m}}$ 
denotes the dual map of $\mathtt{m}\in\mathcal{G}$ from Proposition~\ref{P:dualMapL} and $(r_\mathtt{m})_{\mathtt{m}\in\mathcal{G}}$ are non-negative rates. We assume that $G$ satisfies the summability condition
\begin{align}
\label{sumcond}
\sup_{i\in\Lambda}\sum_{\substack{\mathtt{m}\in\mathcal{G}\\ \mathcal{D}(\mathtt{m})\ni i}}r_\mathtt{m}\big(\vert\mathcal{R}(\mathtt{m}[i])\vert+1\big)<\infty.
\end{align}
Under this condition we can almost surely construct a unique interacting particle system $X=(X_t)_{t\geq0}$ with generator $G$ on $S^\Lambda$ (see \cite[Theorem~4.30]{Swa22}). It turns out (see Theorem~\ref{T:JansBIGtheo} below) that this condition moreover already implies that there exists a non-explosive Markov chain $(Y_t)_{t\geq0}$ with generator $\widehat{G}$ on the countable state space $R^\Lambda_\text{fin}$. We want to prove the following generalisation of Proposition~\ref{P:dual}.
\begin{theorem}[Pathwise monoid duality]
\label{T:dual}
Let $S,R$ and $T$ be finite commutative monoids so that $S$ is $T$-dual to $R$ with duality function $\psi$. Let $G$ and $\widehat{G}$ be the generators from (\ref{monG}) and (\ref{monGhat}) defined via $\mathcal{G}$, a countable collection of maps in $\mathcal{H}_\text{\emph{loc}}(S^\Lambda,S^\Lambda)$ and their unique dual maps from Proposition~\ref{P:dualMapL}. Assume that $G$ satisfies (\ref{sumcond}). 
Fixing a $T\geq0$, we can almost surely construct $X=(X_t)_{t\geq0}$, the process with generator $G$, and $Y=(Y_t)_{t\geq0}$, the process with generator $\widehat{G}$, in such a way that for every $t\in[0,T]$ the random variables $X_t$ and $Y_{T-t}$ are independent and
\begin{align}
\label{pathdual:gen}
[0,T]\ni t\mapsto\Psi\big(X_t,Y_{T-t}^-\big)
\end{align}
is constant, where $Y^-=(Y^-_t)_{t\geq0}$ is the c\`agl\`ad modification of $Y$.
\end{theorem}
By definition, we say that $X$ and $Y$ are \emph{pathwise dual} if they can be constructed in such a way that (\ref{pathdual:gen}) is satisfied. To prove the above result we cite general theory from \cite{Swa22}. 

Let $L$ and $V$ be arbitrary finite sets and let $\mathcal{Y}$ be an arbitrary countable set. As always, we equip $L^\Lambda$ with the product topology and $\mathcal{Y}$ with the discrete one. Let $\varphi:L^\Lambda\times\mathcal{Y}\to V$ be a function. Let $\mathcal{H}$ be a countable collection of local maps in $\mathtt{m}:L^\Lambda\to L^\Lambda$ and assume that every $m\in\mathcal{H}$ has a unique dual map $\hat{\mathtt{m}}:\mathcal{Y}\to\mathcal{Y}$ with respect to $\varphi$. Let $(r_\mathtt{m})_{\mathtt{m}\in\mathcal{H}}$ be non-negative rates and define formal generators $H$ and $\widehat{H}$ in parallel to (\ref{monG}) and (\ref{monGhat}) with $\mathcal{G}$ replaced by $\mathcal{H}$. 
Let $\omega$ denote a Poisson point set on $\mathcal{H}\times\mathbb{R}$ with intensity measure $\rho(\{m\}\times A):=r_m\ell(A)$ $(m\in\mathcal{H},\ A\in\mathcal{B}(\mathbb{R}))$, where $\ell$ denotes the Lebesgue measure. Under condition (\ref{sumcond}), \cite[Theorem~6.16]{Swa22} says that we can almost surely define stochastic flows\footnote{By definition, $(\mathbf{Z}_{s,u})_{z\leq u}$ is a stochastic flow if $\mathbf{Z}_{s,s}$ is the identity map for all $s\in\mathbb{R}$ and if $\mathbf{Z}_{t,u}\circ\mathbf{Z}_{s,t}=\mathbf{Z}_{s,u}$ ($s\leq t\leq u$).} $(\mathbf{X}^+_{s,u})_{s\leq u}$ and $(\mathbf{X}^-_{s,u})_{s\leq u}$ of random continuous maps from $L^\Lambda$ to itself so that, for $s\leq u$, $\omega^+_{s,u}:=\{(m,t)\in\omega:t\in(s,u]\}$ and $\omega^-_{s,u}:=\{(m,t)\in\omega:t\in[s,u)\}$,
\begin{align}
\label{pwlim}
\mathbf{X}^\pm_{s,u}(x)=\lim_{\omega_n\uparrow\omega^\pm_{s,u}}\mathbf{X}^{\omega_n}_{s,u}(x)\qquad(x\in L^\Lambda)
\end{align}
pointwise, where $(\omega_n)_n$ is an arbitrary increasing sequence of finite subsets of $\omega^\pm_{s,u}$ whose union is $\omega^\pm_{s,u}$, and $\mathbf{X}^{\omega_n}_{s,u}$ is the concatenation of all maps in in $\omega_n$ (ordered by the time coordinate $t$).

Let $\widehat{\mathcal{H}}:=\{\hat{m}:m\in\mathcal{H}\}$ and let $\hat{\omega}$ be defined by
\begin{align*}
\hat{\omega}:=\{(\hat{m},-t):(m,t)\in\omega\}.
\end{align*}
Then $\hat{\omega}$ is a Poisson point set on $\widehat{\mathcal{H}}\times\mathbb{R}$ with intensity measure $\hat{\rho}(\{\hat{m}\}\times A):=r_m\ell(A)$ and analogously to above we can almost surely define stochastic flows $(\mathbf{Y}^+_{s,u})_{s\leq u}$ and $(\mathbf{Y}^-_{s,u})_{s\leq u}$ of random continuous maps from $\mathcal{Y}$ to itself so that, for $s\leq u$, $\mathbf{Y}^+_{s,u}$ and $\mathbf{Y}^-_{s,u}$ correspond to pointwise limits as in (\ref{pwlim}), replacing $\omega^+_{s,u}$ by $\hat{\omega}^+_{s,u}=\{(\hat{m},t):(m,t)\in\omega^-_{-u,-s}\}$ and $\omega^-_{s,u}$ by $\hat{\omega}^-_{s,u}=\{(\hat{m},t):(m,t)\in\omega^+_{-u,-s}\}$. The next statement follows from \cite[Theorem~6.20]{Swa22}.
\begin{theorem}[Pathwise dual of an IPS]
\label{T:JansBIGtheo}
Assume that the function $\varphi:L^\Lambda\times\mathcal{Y}\to V$ is continuous if we equip $L^\Lambda\times\mathcal{Y}$ with the product topology, and that it satisfies property (iii) of the definition of duality, i.e.\ that $\varphi(x,y_1)=\varphi(x,y_2)$ for all $x\in L^\Lambda$ implies $y_1=y_2$ ($y_1,y_2\in\mathcal{Y}$). Further assume that $H$ satisfies (\ref{sumcond}). Then there exists a continuous-time Markov chain with generator $\widehat{H}$ that is non-explosive. Moreover, constructing $(\mathbf{X}^\pm_{s,u})_{s\leq u}$ and $(\mathbf{Y}^\pm_{s,u})_{s\leq u}$ as above,
\begin{align}
\label{pwDUALeq}
\varphi\big(\mathbf{X}^\pm_{s,u}(x),y\big)=\varphi\big(x,\mathbf{Y}^\mp_{-u,-s}(y)\big)
\end{align}
holds almost surely simultaneously for all $s\leq u$, $x\in L^\Lambda$ and $y\in\mathcal{Y}$.
\end{theorem}
If two stochastic flows satisfy (\ref{pwDUALeq}) for all $s\leq u$ and for all $x$ and $y$, we say that they are \emph{dual}. Theorem~\ref{T:dual} follows now almost directly from Theorem~\ref{T:JansBIGtheo}.
\begin{proof}[Proof of Theorem~\ref{T:dual}]
First note that Proposition~\ref{P:topdualSR} and the definition of the product topology imply that, by property (ii) of the definition of duality, $\Psi$ from (\ref{Psi}) is also continuous as a function from $S^\Lambda\times R^\Lambda_\text{fin}$ to $T$. Proposition~\ref{P:dualMapL} and Theorem~\ref{T:JansBIGtheo} then show that we can, almost surely, construct stochastic flows $(\mathbf{X}^\pm_{s,u})_{s\leq u}$ and $(\mathbf{Y}^\pm_{s,u})_{s\leq u}$ corresponding to the maps in $\mathcal{G}$ as in Theorem~\ref{T:JansBIGtheo}.

Fix now $T\geq0$ and choose a random variable $X_0$ on $S^\Lambda$ and a random variable $Y_0$ on $R^\Lambda_\text{fin}$, both independent of $(\mathbf{X}^+_{s,u})_{s\leq u}$ and $(\mathbf{Y}^-_{s,u})_{s\leq u}$. Setting
\begin{align*}
X_t:=\mathbf{X}^+_{0,t}(X_0)\quad\text{and}\quad Y_t:=\mathbf{Y}^+_{-T,t-T}(Y_0)\qquad(t\geq0)
\end{align*}
yields by \cite[Proposition~2.9 \& Theorem~4.20]{Swa22} and Theorem~\ref{T:JansBIGtheo} a Markov process $X=(X_t)_{t\geq0}$ with  generator $G$ and a non-explosive continuous-time Markov chain $Y=(Y_t)_{t\geq0}$ with generator $\widehat{G}$. By the construction in \cite[Section~6.4]{Swa22} defining $Y_t^-:=\mathbf{Y}^-_{-T,t-T}(Y_0)$ for $t\geq0$ gives the c\`agl\`ad modification $Y^-:=(Y^-_t)_{t\geq0}$ of $Y$. Using the duality of the stochastic flows, i.e.\ (\ref{pwDUALeq}), one then has for all $s,u\in\mathbb{R}$ satisfying $0\leq s\leq u\leq T$ that
\begin{align*}
\begin{split}
\Psi\big(X_s,Y_{T-s}^-\big)&=\Psi\big(\mathbf{X}^+_{0,s}(X_0),\mathbf{Y}^-_{-T,-s}(Y_0)\big)=\Psi\big(\mathbf{X}^+_{0,s}(X_0),\mathbf{Y}^-_{-u,-s}\circ\mathbf{Y}^-_{-T,-u}(Y_0)\big)\\
&=\Psi\big(\mathbf{X}^+_{s,u}\circ\mathbf{X}^+_{0,s}(X_0),\mathbf{Y}^-_{-T,-u}(Y_0)\big)=\Psi\big(\mathbf{X}^+_{0,u}(X_0),\mathbf{Y}^-_{-T,-u}(Y_0)\big)\\
&=\Psi\big(X_u,Y_{T-u}^-\big),
\end{split}
\end{align*}
i.e.\ the function in (\ref{pathdual:gen}) is constant, and the proof is complete.
\end{proof}

Applying the general theory to the 2CP we prove Proposition~\ref{P:dual}.
\begin{proof}[Proof of Proposition~\ref{P:dual}]
As already mentioned in Section~\ref{S:Int}, $U=(U,\veebar)$ is indeed a monoid. Next one computes $\mathcal{H}(U,U)$ and $\mathcal{H}(U,M)$, with $M:=(\{-1,0,1\},\,\cdot\,)$, where $\cdot$ denotes the usual multiplication in $\mathbb{R}$. To compute $\mathcal{H}(U,U)=\{(o,o),(o,\text{id}),(\text{id},o),(\text{id},\text{id})\}$ one can apply \cite[Lemma~6]{LS22}, noting that $U=M_1\times M_2$, where $M_1:=(\{0,1\},\vee)$ and $M_2:=(\{0,1\},\oplus)$, and checking that $\Hi(M_i,M_j)=\{o,\text{id}\}$ if $i=j$ and $=\{o\}$ if $i\neq j$ $(i,j\in\{1,2\})$. To compute $\mathcal{H}(U,M)$ one can apply the same result, computing first $\mathcal{H}(M_1,M)=\{1,\gamma_1\}$ and $\mathcal{H}(M_2,M)=\{1,\gamma_2\}$, where $1$ is the function constantly 1, and $\gamma_1$ and $\gamma_2$ are the functions from (\ref{gammas}). Using the definition of duality one then confirms that $U$ is $M$-dual to itself with respect to $\psi$ from (\ref{psidefloc}).

Having computed $\mathcal{H}(U,U)$ one directly concludes that all its maps are self-dual as $o$ and $\text{id}$ are always self-dual. All maps in (\ref{prodmaps}) (that are used in the definition of $G_\veebar$ in (\ref{2G})) can be written as in (\ref{DefMATL_H}) with $\Delta$ from (\ref{EnotID}) finite, so Lemma~\ref{L:L_HwithMAT} implies that they are elements of $\mathcal{H}_\text{loc}(U^\Lambda,U^\Lambda)$, with $\Lambda=\mathbb{Z}^d$ and $U^\Lambda$ being, as always, equipped with the product topology. Proposition~\ref{P:dualMapL} shows that $G_\veebar$ can play the role of both $G$ and $\widehat{G}$ from (\ref{monG}) and (\ref{monGhat}). One quickly verifies that (\ref{sumcond}) holds and the claim follows from Theorem~\ref{T:dual}.
\end{proof}

One can check that the monoid $U$ is isomorphic to $M_{23}$ from \cite[Appendix~A.1]{LS22} and the monoid $M=(\{-1,0,1\},\,\cdot\,)$ is isomorphic to $M_5$ from \cite[Section~5.1]{LS22}. The function $\psi$ is denoted in \cite{LS22} as $\psi_{235}$ and the fact that $U$ is $M$-dual to itself can be found in the table in \cite[Appendix~A.2]{LS22}. The fact that $\mathcal{H}(U,U)=\{(o,o),(o,\text{id}),(\text{id},o),(\text{id},\text{id})\}$ is, by \cite[Proposition~4]{LS22}, encoded in the duality function $\psi_{23}$ from \cite[Appendix~A.2]{LS22}.

%
%

\section{Informativeness and representations}
\label{S:Info}

In this subsection Lemma~\ref{L:info} is proved. In fact, as already stated in the outline, we are going to prove a more general result and we are going to investigate the open task to classify the monoid dualities from \cite{LS22} that determine the law of processes uniquely. Let, as in the section above, $(S,\odot)$, $(R,\raisebox{-0.35mm}{$\boxdot$})$ and $(T,\otimes)$ be commutative finite monoids and assume that $S$ is $T$-dual to $R$ with duality function $\psi:S\times R\to T$. Let $\Lambda$ be countable, let $\mathbb{V}$ be a finite dimensional real or complex vector space and let $V$ be an arbitrary measurable space. 

Towards the goal of classification we give the following definitions. For an arbitrary index set $I$ we call a family $(f_i)_{i\in I}$ of measurable functions $f_i:S^\Lambda\to\mathbb{V}$ \emph{distribution determining} if, for two random variables $X$ and $X^\prime$ on $S^\Lambda$,
\begin{align*}
\mathbb{E}[f_i(X)]=\mathbb{E}[f_i(X^\prime)]\quad\forall i\in I\qquad&\text{implies}\qquad X\stackrel{d}{=}X^\prime,
\intertext{where $\stackrel{d}{=}$ denotes equality in distribution. Similarly, we call a family $(g_i)_{i\in I}$ of measurable functions $g_i:S^\Lambda\to V$ \emph{weakly distribution determining} if}
g_i(X)\stackrel{d}{=}g_i(X^\prime)\quad\forall i\in I\qquad&\text{implies}\qquad X\stackrel{d}{=}X^\prime.
\end{align*}
The first of the two definition is already widely used (compare \cite{Swa22}), while the second one we introduce here newly. 

A family $(f_i)_{i\in I}$ of functions $f_i:S^\Lambda\to\mathbb{V}$ that is distribution determining is clearly also weakly distribution determining. The reverse implication is not true in general, but holds in the following special case. Recall that $v_1,\ldots,v_n\in\mathbb{V}$ are called \emph{affinely independent} if 
\begin{align*}
\sum_{k=1}^n \lambda_k v_k=0\ \text{with scalars}\ \lambda_1,\ldots,\lambda_n\ \text{s.t.}\ \sum_{k=1}^n\lambda_k=0\quad\text{implies}\quad\lambda_1=\ldots=\lambda_n=0.
\end{align*}
\begin{proposition}[Equality of notions]
\label{P:affindep}
Let $(f_i)_{i\in I}$ be a family of functions $f_i:S^\Lambda\to\{v_1,\ldots,v_n\}\subset\mathbb{V}$. If $v_1,\ldots,v_n$ are affinely independent, then $(f_i)_{i\in I}$ is distribution determining if and only if it is weakly distribution determining.
\end{proposition}
\begin{proof}
Comparing the definitions it suffices to show for fixed $i\in I$ that, under the assumption of the proposition, $\mathbb{E}[f_i(X)]=\mathbb{E}[f_i(X^\prime)]$ implies $f_i(X)\stackrel{d}{=}f_i(X^\prime)$. As the set $\{v_1,\ldots,v_n\}$ is finite, the condition $\mathbb{E}[f_i(X)]=\mathbb{E}[f_i(X^\prime)]$ is equivalent to writing
\begin{align*}
\sum_{k=1}^n v_k\big(\mathbb{P}[f_i(X)=v_k]-\mathbb{P}[f_i(X^\prime)=v_k]\big)=0.
\end{align*}
But as $v_1,\ldots,v_n$ are affinely independent, then also
\begin{align*}
\mathbb{P}[f_i(X)=v_k]-\mathbb{P}[f_i(X^\prime)=v_k]=0\qquad(k=1,\ldots,n),
\end{align*}
i.e.\ $f_i(X)$ and $f_i(X^\prime)$ are equal in distribution.
\end{proof}

Let now $\Psi:S^\Lambda\times R^\Lambda_\text{fin}\to T$ be the function from (\ref{Psi}). In parallel to \cite{Swa22} we say that $\Psi$ is \emph{weakly informative} if
\begin{align}
\label{INFPfcts}
\big(\Psi(\,\cdot\,,\mathbf{y})\big)_{\mathbf{y}\in R^\Lambda_\text{fin}}
\end{align}
is weakly distribution determining. If the monoid $T$ is also a subset of a real or complex vector space, we say that $\Psi$ is \emph{informative} if the functions in (\ref{INFPfcts}) are distribution determining. We prove the following result.
\begin{proposition}[Informativeness of $\Psi$]
\label{P:genINFO}
Under the assumptions of this subsection $\Psi$ is informative if $T$ is a sub-monoid of $(\mathbb{C},\,\cdot\,)$, where $\cdot$ denotes the usual multiplication.
\end{proposition}
It is easy to see that all finite sub-monoids of $(\mathbb{C},\,\cdot\,)$ (apart from $(\{0\},\,\cdot\,)$) consist of the multiplicative group of $n$-th roots of unity for some $n\in\mathbb{N}$, either with or without an added 0. Those with cardinality up to four are named $M_0,M_1, M_2, M_5, M_7, M_{18}$ and $M_{26}$ in our paper \cite{LS22}, so by Proposition~\ref{P:genINFO} all duality functions from \cite{LS22} that take values in these monoids are informative. In particular, setting $(T,\otimes)=(\{1,-1,0\},\,\cdot\,)$ and  $(S,\odot)=(R,\raisebox{-0.35mm}{$\boxdot$})=(U,\veebar)$, Proposition~\ref{P:genINFO} implies Lemma~\ref{L:info}. 

To prove Proposition~\ref{P:genINFO} we use a Stone-Weierstrass argument. Let $\mathcal{C}(\mathcal{X},\mathcal{Y})$ denote the space of continuous functions from space $\mathcal{X}$ to space $\mathcal{Y}$. We say that $\mathcal{H}\subset\mathcal{C}(\mathcal{X},\mathcal{Y})$ \emph{separates points} if for $x,x^\prime\in\mathcal{X}$ with $x\neq x^\prime$ there exists $f\in\mathcal{H}$ such that $f(x)\neq f(x^\prime)$. Moreover, we say that $\mathcal{G}\subset\mathcal{C}(\mathcal{X},\mathbb{C})$ is \emph{self-adjoint} if $f\in\mathcal{G}$ implies $\overline{f}\in\mathcal{G}$, where $\overline{f}(x):=\overline{f(x)}$ $(x\in\mathcal{X})$, the complex conjugate of $f(x)$.
\begin{lemma}[Application of Stone-Weierstrass]
\label{L:Stone-Weier}
Let $E$ be a compact metrizable space. Assume that $\mathcal{G}\subset\mathcal{C}(E,\mathbb{C})$ separates points and is closed under products. Then $\mathcal{G}$ is distribution determining.
\end{lemma}
\begin{proof}
The statement with $\mathbb{C}$ replaced by $\mathbb{R}$ is proved in \cite[Lemma~4.37]{Swa22}. Note that
\begin{align}
\label{compconj}
\mathbb{E}\big[f(X)\big]=\mathbb{E}\big[f(X^\prime)\big]\quad\text{implies}\quad\mathbb{E}\left[\overline{f}(X)\right]=\mathbb{E}\left[\overline{f}(X^\prime)\right]\qquad(f\in\mathcal{G}),
\end{align}
as $\mathbb{E}\big[\overline{f}(X)\big]=\overline{\mathbb{E}\big[f(X)\big]}$, where $X$ and $X^\prime$ are random variables on $E$. We can enlarge $\mathcal{G}$ with the constant function 1, take linear combinations and convex conjugates and receive an algebra $\mathcal{H}\supset\mathcal{G}$ that is closed under products, self-adjoint and separates points. If $\mathbb{E}\big[f(X)\big]=\mathbb{E}[f(X^\prime)]$ for all $f\in\mathcal{G}$ then also $\mathbb{E}\big[f(X)\big]=\mathbb{E}[f(X^\prime)]$ for all $f\in\mathcal{H}$ by the linearity of the integral and (\ref{compconj}). We then can apply the complex version of the Stone-Weierstrass theorem and continue as in the proof of \cite[Lemma~4.37]{Swa22}.
\end{proof}

\begin{proof}[Proof of Proposition~\ref{P:genINFO}]
By definition, we have to prove that the family
\begin{align*}
\mathcal{G}:=\big(\Psi(\,\cdot\,,\mathbf{y})\big)_{\mathbf{y}\in R^\Lambda_\text{fin}}
\end{align*}
is distribution determining.

By Tychonoff's theorem, the space $S^\Lambda$, equipped with the product topology, is a compact metrizable space. The fact that $\mathcal{G}$ is closed under products follows from the duality between $S^\Lambda$ and $R^\Lambda_\text{fin}$: Property (i) in the definition of duality implies that
\begin{align*}
\Psi(\mathbf{x},\mathbf{y}_1)\Psi(\mathbf{x},\mathbf{y}_2)=\Psi(\mathbf{x},\mathbf{y}_1\hspace{0.2em}\raisebox{-0.35mm}{$\boxdot$}\hspace{0.2em}\mathbf{y}_2)\qquad(\mathbf{x}\in S^\Lambda,\ \mathbf{y}_1,\mathbf{y}_2\in R^\Lambda_\text{fin}).
\end{align*}
\noi
The fact that $\mathcal{G}$ separates points follows directly from property (ii) of the definition of (topological) duality. Applying Lemma~\ref{L:Stone-Weier} then yields Proposition~\ref{P:genINFO}.
\end{proof}

To further investigate the case in which the monoid $T$ can not naturally be written as a sub-monoid of $(\mathbb{C},\,\cdot\,)$, we provide some additional notions. The reader that is just concerned with the 2CP may skip ahead to the next section.

A \emph{multiplicative representation} of a commutative monoid $(M,+)$ with neutral element 0 is a map $\ga:M\to\A$, where $(\A,+,\,\cdot\,)$ is a unital commutative algebra with unit $I$, so that $\ga(x+y)=\ga(x)\cdot\ga(y)$ and $\ga(0)=I$. Then $\gamma(M)=\{\ga(x):x\in M\}$ is a sub-monoid of $(\A,\,\cdot\,)$ and $\ga:M\to\gamma(M)$ is a homomorphism. We say that $\ga$ is \emph{faithful} if this is an isomorphism.

We again consider the function $\Psi:S^\Lambda\times R^\Lambda_\text{fin}\to T$ from (\ref{Psi}). Recall that, by Proposition~\ref{P:topdualSR}, under the usual assumptions on $S,R$ and $T$ (see the beginning of this section), $S^\Lambda$ is $T$-dual to $R^\Lambda_\text{fin}$ with duality function $\Psi$. If now $\ga:T\to\A$ is a faithful multiplicative representation, we equip the finite monoids $T$ and $\gamma(T)$ with the discrete topology and it follows from the definition of duality that $S^\Lambda$ is also $\gamma(T)$-dual to $R^\Lambda_\text{fin}$ with duality function $\gamma\circ\Psi$. If $\gamma\circ\Psi$ is informative we say that $\gamma\circ\Psi$ is a \emph{good multiplicative representation} of $\Psi$. Proposition~\ref{P:affindep} and the faithfulness of $\gamma$ imply that $\gamma\circ\Psi$ is a good representation of $\Psi$ if $\Psi$ is weakly informative as long as the elements of $\gamma(T)$ are affinely independent. The next result states that we can always find such a good multiplicative representation of a weakly informative duality function, so weak informativeness is basically all we need in practice.

\begin{proposition}[Existence of good representations]
\label{P:EgoodREP}
\hspace{-0.5em} Under the assumptions of this subsection there exist a finite dimensional real unital commutative algebra $\A$ and a faithful representation $\ga:T\to\A$ such that $\gamma\circ\Psi$ is informative if $\Psi$ is weakly informative.
\end{proposition}
\begin{proof}
Let $\mathbb{R}^T$ be the space of all functions mapping from $T$ to $\mathbb{R}$. The space $(\mathbb{R}^T,+)$, where + denotes the usual (pointwise) sum of real-valued functions, is a finite dimensional real vector space on which we can define the product $\ast$ as
\begin{align*}
(g\ast h)(x):=\sum_{y,z\in T}g(y)h(z)\mathbbm{1}_{\{x\}}(y\otimes z)\qquad\big(g,h\in\mathbb{R}^T,\ x\in T\big),
\end{align*}
where the sum is the usual sum in $\mathbb{R}$ and $\mathbbm{1}$ denotes the indicator function. One readily checks that this makes $(\mathbb{R}^T,+,\ast)$ a finite dimensional real unital algebra with unit $\mathbbm{1}_{\{0\}}$. Defining $\gamma:T\to\mathbb{R}^T$ as $\gamma(x)=\mathbbm{1}_{\{x\}}$ ($x\in T$) then gives a faithful multiplicative representation of $T$ and clearly the elements of $\gamma(T)$ are affinely independent. The claim then follows from Proposition~\ref{P:affindep} and the faithfulness of $\gamma$ as stated above.
\end{proof}

By the above proposition we can reformulate the classification problem by asking to classify general duality functions (that do not map into sub-monoids of $(\mathbb{C},\,\cdot\,)$) into the classes ``weak informative" and ``not weak informative". This remains an open problem.

We end this section with an additional observation. While $\R^T$ from the proof of Proposition~\ref{P:EgoodREP} is a $\vert T\vert$-dimensional vector space, Proposition~\ref{P:genINFO} implies that for large $T$ also representations in lower dimensional spaces can be good, even if the elements of $\gamma(T)$ are not affinely independent. As it is in practice often easier to work in a lower dimensional space, there can exist ``better" representations of weakly informative duality functions than the one from Proposition~\ref{P:EgoodREP}. In light of Proposition~\ref{P:genINFO} one might even hope that $\gamma\circ\Psi$ is always a good representation of a weakly informative $\Psi$ as long as $\gamma$ is faithful. This, however, is not true and we provide a counterexample below.

We again consider the monoid $(U,\veebar)$ defined in Section~\ref{S:PoI}. 
From \cite[Appendix~A.2]{LS22} we know that there also exists the ``local" duality function $\psi_{23}$ mapping from $U\times U$ back into $U$. Reordering the elements of $M_{23}$ as in the present paper (i.e.\ as in $U$) one has that
\begin{align*}
\psi_{23}(x,y)=\big(\psi_1(x_1,y_1),\psi_2(x_2,y_2)\big)\qquad\big(x=(x_1,x_2),y=(y_1,y_2)\in U\big),
\end{align*}
where $\psi_1$ and $\psi_2$ are the ``local" additive and cancellative duality function, defined in (\ref{psi1&2}). It follows from (\ref{frompsi1&2topsi}) that
\begin{align}
\label{23impliespsi}
\psi_{23}(x,y)=\psi_{23}(v,w)\qquad\text{implies}\qquad\psi(x,y)=\psi(v,w)\qquad\big(x,y,v,w\in U\big).
\end{align}
We define a ``global" duality function $\boldsymbol{\psi}_{23}:\mathcal{U}\times\mathcal{U}_\text{fin}\to U$ as in (\ref{psidef}), but for $\psi_{23}$ instead of $\psi$ and with the ``product" taken in $U$. It follows from (\ref{23impliespsi}) that for two random variables $X,X^\prime$ on $\mathcal{U}$ and for $y\in\mathcal{U}_\text{fin}$,
\begin{align*}
\boldsymbol{\psi}_{23}(X,y)\stackrel{d}{=}\boldsymbol{\psi}_{23}(X^\prime,y)\quad\text{implies}\quad\boldsymbol{\psi}(X,y)\stackrel{d}{=}\boldsymbol{\psi}(X^\prime,y),
\end{align*}
and, due to the informativeness of $\boldsymbol{\psi}$, the duality function $\boldsymbol{\psi}_{23}$ is weakly informative. Defining now $\gamma:U\to\R^2$ as
\begin{align*}
\gamma(x):=\big(\gamma_1(x_1),\gamma_2(x_2)\big)\qquad(x=(x_1,x_2)\in U),
\end{align*}
with $\gamma_1,\gamma_2$ defined in (\ref{gammas}), yields a faithful multiplicative representation of $U$ in $\R^2$, viewed as a unital algebra equipped with pointwise multiplication.

However, $\gamma\circ\psib_{23}$ is \emph{not} a good representation of $\psib_{23}$. For example, the random variables $X,X^\prime$ on $\mathcal{U}$ with
\begin{align*}
&\mathbb{P}[X(i)=(0,0)]=\mathbb{P}[X^\prime(i)=(0,0)]=1\ \text{for}\ i\in\mathbb{Z}^d\setminus\{0\},\\[0.5em]
&\mathbb{P}[X(0)=x]=\frac{1}{4}\ \text{for all}\ x\in U\qquad \mathbb{P}[X^\prime(0)=x]=\begin{cases}
\frac{1}{2}&\text{if}\ x\in\{(0,0),(1,1)\},\\
0&\text{else},
\end{cases}
\end{align*}
show that $\gamma\circ\boldsymbol{\psi}_{23}:\mathcal{U}\times\mathcal{U}_\text{fin}\to\mathbb{R}^2$ is not informative. Here $0\in\mathbb{Z}^d$ denotes the origin.


%
%

\section{The main convergence result}
\label{S:weak}

In this section we prove Theorem~\ref{T:invlaw}. Recall that $\mathcal{T}$ denotes the space of all functions $z:\mathbb{Z}^d\to T=\{0,1\}$ and recall the definition of $\mathcal{T}_\text{fin}$ in (\ref{T_fin}). For $z\in\mathcal{T}$ we shortly write $\vert z\vert:=\vert\{i\in\Z^d:z(i)=1\}\vert$. We are going to use several auxiliary lemmas to prove Theorem~\ref{T:invlaw}. The first one is \cite[Lemma~6.37]{Swa22}. 
The symbol $\wedge$ denotes the pointwise minimum, i.e.\ $(z_1\wedge z_2)(i)=\min\{z_1(i),z_2(i)\}$ for $i\in\mathbb{Z}^d,\ z_1,z_2\in\mathcal{T}$.
\begin{lemma}[Non-zero intersection: CP]
\label{II:L2}
Let $Z=(Z_t)_{t\geq0}$ be a \emph{CP(}$\lambda,\delta$\emph{)} $(\lambda>0,\ \delta\geq0)$ with non-trivial shift-invariant initial distribution. Given $\varepsilon>0$, for each time $s>0$ there exists an $N_\text{\emph{CP}}\in\mathbb{N}$ such that for any $z\in\mathcal{T}$ with $\vert z\vert\geq N_\text{\emph{CP}}$ one has
\begin{align*}
\mathbb{P}\big(Z_s\wedge z=\underline{0}\big)\leq\varepsilon.
\end{align*}
\end{lemma}
Additionally we are going to use the following application of \cite[Corollary~9]{SS08a}. As \cite[Corollary~9]{SS08a} is not stated in the most accessible form we devote Appendix~\ref{S:App} to showing how the result below follows from it. Instead of using the result below we could have also followed the strategy of the proof of \cite[Theorem~1.2]{BDD91}. There the authors use the graphical representation of the cCP explicitly to work around the statement below.
\begin{lemma}[Parity indeterminacy: cCP]
Let\label{L:indeter} $Z=(Z_t)_{t\geq0}$ be a \emph{cCP}$(\la,\de)$ $(\la>0,\ \delta\geq0)$ with non-trivial shift-invariant initial distribution. Given $\eps>0$, for each time $s>0$ there exists an $N_\text{\emph{cCP}}\in\mathbb{N}$ such that for any $z\in\mathcal{T}_\text{\emph{fin}}$ with $\vert z\vert\geq N_\text{\emph{cCP}}$ one has
\begin{align*}
\left|\P\big[|Z_s\wedge z\vert\ \text{is odd}\,\big]-\frac{1}{2}\right|\leq\eps.
\end{align*}
\end{lemma} 
Finally, the following result extends \cite[Lemma~6.36]{Swa22} and \cite[Lemma~2.1]{BDD91}.
\begin{lemma}[Extinction or unbounded growth]
\label{L:0vsINFTY}
Let $Z=(Z_t)_{t\geq0}$ be either a \emph{CP(}$\lambda,\delta$\emph{)} or a \emph{cCP(}$\lambda,\delta$\emph{)} $(\lambda,\delta\geq0,\ \lambda+\delta>0)$. For each $z\in\mathcal{T}_\text{\emph{fin}}$ and $N\in\mathbb{N}$ one has  
\begin{align}
\label{0vsINFTYweak}
\lim_{t\to\infty}\mathbb{P}^z[0<\vert Z_t\vert<N]=0.
\end{align}
\end{lemma}
\begin{proof}
If $z=\underline{0}$ the statement is trivial, so let $z\in\mathcal{T}_\text{fin}\setminus\{\underline{0}\}$. In the case $\lambda,\delta>0$ \cite[Lemma~6.36]{Swa22} and \cite[Lemma~2.1]{BDD91} imply
\begin{align}
\label{0vsINFTYstrong}
\mathbb{P}^z\big[\exists t\geq0:Z_t=\underline{0}\ \text{or}\ \vert Z_t\vert\to\infty\ \text{as}\ t\to\infty\big]=1
\end{align}
for the CP and the cCP, respectively, and (\ref{0vsINFTYstrong}) clearly implies (\ref{0vsINFTYweak}). In fact, the two proofs are just reformulations of each other, both based on L\'evy's 0-1 law. 

In the case $\lambda=0,\ \delta>0$ there is no difference between a CP and a cCP and
\begin{align*}
\mathbb{P}^z\big[\exists t\geq 0:Z_t=\underline{0}\big]=\lim_{t\to\infty}\mathbb{P}^z\big[Z_t=\underline{0}\big]=\lim_{t\to\infty}\big(1-e^{-\delta t}\big)^{\vert z\vert}=1
\end{align*}
since $\underline{0}$ is absorbing. This implies (\ref{0vsINFTYstrong}) and hence also (\ref{0vsINFTYweak}).

In the case $\lambda>0,\ \delta=0$, and if $Z$ is a CP, the function $t\mapsto\vert Z_t\vert$ is non-decreasing, hence it converges in $\mathbb{N}\cup\{\infty\}$. Let $N\in\mathbb{N}$. One has
\begin{align}
\label{CParg}
\mathbb{P}^z\big[\lim\nolimits_{t\to\infty}\vert Z_t\vert\leq N\big]=1-\mathbb{P}^z\big[\exists t\geq 0:\vert Z_t\vert>N\big]=1-\lim_{t\to\infty}\mathbb{P}^z[\vert Z_t\vert>N]=0
\end{align}
as choosing a suitable sequence of neighbours and neighbours of neighbours of the infected individuals in $z$ yields that
\begin{align*}
\mathbb{P}^z[\vert Z_t\vert>N]\geq\left(1-\mathbbm{1}_{\{\vert z\vert\leq N\}}e^{-\frac{\lambda t}{N+1-\vert z\vert}}\right)^{N+1-\vert z\vert}
\end{align*}
for $t>0$. Here, in the case that $\vert z\vert\leq N$, we have divided time into $N+1-\vert z\vert$ subintervals and used the fact that $1-e^{-\lambda t}$ is the probability to infect a previously chosen neighbour of an infected individual during a time interval of length $t$. Finally, (\ref{CParg}) implies that
\begin{align*}
\mathbb{P}^z\big[\vert Z_t\vert\to\infty\ \text{as}\ t\to\infty\big]&=1-\mathbb{P}^z\big[\exists N\in\mathbb{N}:\lim\nolimits_{t\to\infty}\vert Z_t\vert=N\big]\\
&\geq1-\sum_{N\in\mathbb{N}}\mathbb{P}^z\big[\lim\nolimits_{t\to\infty}\vert Z_t\vert\leq N\big]=1,
\end{align*}
again implying (\ref{0vsINFTYstrong}) and hence also (\ref{0vsINFTYweak}).

To treat the cCP in the case $\lambda>0,\ \delta=0$, we use \cite[Theorem~1.3]{BDD91}. It says that a cCP(1,0), started in \emph{any} initial state other than $\underline{0}$, converges weakly to the product law assigning probability 1/2 to both 0 and 1 at every node. By changing the time scale the same holds for a cCP($\lambda,0$) with an arbitrary $\lambda>0$. Let $N\in\mathbb{N}$ and $\varepsilon>0$. Choose now an $M=M(N,\varepsilon)>N$ so that $p_N:=\mathbb{P}[X\leq N]<\varepsilon$ if $X$ is a binomially distributed random variable with parameters $n=M$ and $p=1/2$. Additionally, choose an arbitrary $x\in\mathcal{T}_\text{fin}$ with $\vert x\vert=M$.  Then, by the weak convergence,
\begin{align*}
\limsup_{t\to\infty}\mathbb{P}^z\big[\vert Z_t\vert\leq N\big]\leq\lim_{t\to\infty}\mathbb{P}^z\big[\vert Z_t\wedge x\vert\leq N\big]=p_N<\varepsilon,
\end{align*}
implying $\lim_{t\to\infty}\mathbb{P}^z\big[\vert Z_t\vert\leq N\big]=0$ (i.e.\ convergence in probability to $\infty$). Thus (\ref{0vsINFTYweak}) holds.
\end{proof}

Using the three lemmas above we are able to prove Theorem~\ref{T:invlaw}. 


\begin{proof}[Proof of Theorem \ref{T:invlaw}]
Let $Y=(Y^1,Y^2)=(Y^1_t,Y^2_t)_{t\geq0}$ be an independent copy of the 2CP $X=(X^1,X^2)=(X^1_t,X^2_t)_{t\geq0}$ in the formulation of the theorem, but started in the deterministic state $y=(y_1,y_2)\in\mathcal{U}_\text{fin}$. Due to the informativeness of $\boldsymbol{\psi}$ and the compactness of $\mathcal{U}$, the set $\mathcal{G}$ from the proof of Lemma~\ref{L:info} is also convergence determining, i.e.\ showing
\begin{align}
\label{keyLIM}
\lim_{t\to\infty}\mathbb{E}\big[\boldsymbol{\psi}(X_t,y)\big]=\mathbb{P}^y\left[\exists t\geq0:Y_t=\underline{(0,0)}\right]
\end{align}
for all $y\in\mathcal{U}_\text{fin}$ implies (\ref{nuconv}) (compare \cite[Lemma~4.38]{Swa22}). If $y=\underline{(0,0)}$, (\ref{keyLIM}) follows trivially from the definition of $\boldsymbol{\psi}$, so assume $y\neq\underline{(0,0)}$. 
We set
\begin{align*}
\lambda_1:=\lambda+\lambda^\vee,\quad\delta_1:=\delta+\delta^\vee,\quad\lambda_2:=\lambda+\lambda^\oplus,\quad\delta_2:=\delta+\delta^\oplus,
\end{align*}
so that $X^1$ and $Y^1$ are both a CP($\lambda_1,\delta_1$), and $X^2$ and $Y^2$ are both a cCP($\lambda_2,\delta_2$). Assume, for now, that $\lambda_1,\lambda_2>0$, so that all three auxiliary lemmas above are applicable. 
Let $\varepsilon>0$ be arbitrary. 
Choose $N_\text{CP}$ and $N_\text{cCP}$ as in Lemma~\ref{II:L2} and Lemma~\ref{L:indeter} in dependence of the chosen $\varepsilon$, $s=1$, and the model parameters. 
Fix $t>0$. We have, using the duality equation (\ref{dual}) and the law of total expectation, that
\begin{align}
\begin{split}
\label{E[PSI]first}
&\mathbb{E}[\boldsymbol{\psi}(X_{t+1},y)]\\
&=\mathbb{E}[\boldsymbol{\psi}(X_{1},Y_{t})]\\
&=\mathbb{E}\big[\boldsymbol{\psi}(X_{1},Y_{t})\mid Y_{t}^1=Y^2_{t}=\un 0\big]\mathbb{P}^y\big[Y_{t}^1=Y^2_{t}=\un 0\big]\\
&\hspace{1em}+\mathbb{E}\big[\boldsymbol{\psi}(X_{1},Y_{t})\mid Y_{t}^1=\un 0,\ 0<\vert Y_{t}^2\vert<N_\text{cCP}\big]\underbrace{\mathbb{P}^y\big[Y_{t}^1=\un 0,\ 0<\vert Y_{t}^2\vert<N_\text{cCP}\big]}_{=:p_1(y,t)}\\
&\hspace{1em}+\underbrace{\mathbb{E}\big[\boldsymbol{\psi}(X_{1},Y_{t})\mid Y^1_{t}=\un 0,\ \vert Y^2_{t}\vert\geq N_\text{cCP}\big]}_{=:E_1(y,t)}\mathbb{P}^y\big[Y^1_{t}=\un 0,\ \vert Y^2_{t}\vert\geq N_\text{cCP}\big]\\
&\hspace{1em}+\mathbb{E}\big[\boldsymbol{\psi}(X_{1},Y_{t})\mid 0<\vert Y^1_{t}\vert<N_\text{CP}\big]\underbrace{\mathbb{P}^y\big[0<\vert Y^1_{t}\vert<N_\text{CP}\big]}_{=:p_2(y,t)}\\
&\hspace{1em}+\underbrace{\mathbb{E}\big[\boldsymbol{\psi}(X_{1},Y_{t})\mid \vert Y^1_{t}\vert\geq N_\text{CP}\big]}_{=:E_2(y,t)}\mathbb{P}^y\big[\vert Y^1_{t}\vert\geq N_\text{CP}\big].
\end{split}
\end{align}

Depending on the choice of the model parameters and $y$, the deterministic initial state of $Y$, it might happen that some of the events on which we condition above have probability zero. The cases that either $y_1=\un 0$ or $y_2=\un 0$, or the monotonely coupled case $\delta_\vee=\lambda_\oplus=0$ when $y$ satisfies $y(i)\neq(0,1)$ for all $i\in\mathbb{Z}^d$ are such examples. In these cases we define the corresponding conditioned expectation (arbitrarily) to equal 1. Due to the zero probability the line in (\ref{E[PSI]first}) where it occurs then drops out, and for the remaining ones we can argue as below.

From the definition of $\psib$ it is clear that $\mathbb{E}\big[\boldsymbol{\psi}(X_{1},Y_{t})\mid Y_{t}^1=Y^2_{t}=\un 0\big]=1$ and
\begin{align*}
\mathbb{P}^y\big[Y_{t}^1=Y^2_{t}=\un 0\big]\nearrow\mathbb{P}^y\left[\exists t\geq0:Y_t=\underline{(0,0)}\right]
\end{align*}
as $t\to\infty$. Moreover, Lemma~\ref{L:0vsINFTY} implies that
\begin{align*}
\lim_{t\to\infty}p_1(y,t)=\lim_{t\to\infty}p_2(y,t)=0.
\end{align*}
As in the proof of \cite[Theorem~6.35]{Swa22} we use Lemma~\ref{II:L2} to compute that
\begin{align}
\begin{split}
\label{useII:L2}
\vert E_2(y,t)\vert&=\left\vert\mathbb{P}\big[\boldsymbol{\psi}(X_{1},Y_{t})=1\mid \vert Y^1_t\vert\geq N_\text{CP}\big]-\mathbb{P}\big[\boldsymbol{\psi}(X_{1},Y_{t})=-1\mid \vert Y^1_t\vert\geq N_\text{CP}\big]\right\vert\\
&\leq\mathbb{P}\big[\boldsymbol{\psi}(X_{1},Y_{t})\neq0\mid \vert Y^1_t\vert\geq N_\text{CP}\big]\\
&=\mathbb{P}\big[X_{1}^1\wedge Y^1_{t}=\underline{0}\mid \vert Y^1_t\vert\geq N_\text{CP}\big]\leq\varepsilon
\end{split}
\end{align}
by the choice of $N_\text{CP}$. For $E_1(y,t)$ one has that
\begin{align*}
E_1(y,t)&=1-2\mathbb{P}\big[\boldsymbol{\psi}(X_{1},Y_{t})=-1\mid Y^1_{t}=\un 0,\ \vert Y^2_{t}\vert\geq N_\text{cCP}\big]\\
&=1-2\mathbb{P}\big[\vert X_{1}^2\wedge Y^2_{t}\vert\ \text{is odd}\ \big\vert\ Y^1_{t}=\underline{0},\ \vert Y^2_{t}\vert\geq N_\text{cCP}\big]
\end{align*}
and, due to the independence of $X$ and $Y$, we can apply Lemma~\ref{L:indeter} and conclude that
\begin{align*}
\vert E_1(y,t)\vert\leq2\varepsilon.
\end{align*}
Plugging then back into (\ref{E[PSI]first}) and computing the limit inferior and the limit superior, one concludes (\ref{keyLIM}) as $\varepsilon$ was arbitrary.

To finish the proof we consider the case that $\lambda_1=0$ and/or $\lambda_2=0$. By assumption, $\lambda_i$ $(i\in\{1,2\})$ can only equal zero if $\delta_i>0$. The idea is to still use (\ref{E[PSI]first}), where we used $\lambda_1>0$ for the treatment of $E_2(y,t)$ and $\lambda_2>0$ for the treatment of $E_1(y,t)$. However, if $\lambda_1=0$, then $Y^1$ is a CP($0,\delta_1$) with $\delta_1>0$, so the number of infected individuals can only decrease. Choosing $N_\text{CP}:=\vert y_1\vert+1$ makes the line in (\ref{E[PSI]first}) in which $E_2(y,t)$ appears vanish. Analogously, choosing $N_\text{cCP}:=\vert y_2\vert+1$ makes the line in which $E_1(y,t)$ appears vanish if $\lambda_2=0$. For the rest of the terms one then can argue as above.

We conclude that in all cases (\ref{keyLIM}) holds, thus also (\ref{nuconv}) as explained above. Lastly, it is well-known (compare \cite[Lemma~4.40]{Swa22}) that (\ref{nuconv}) implies that $\nu$ is indeed invariant and the proof is complete.

\end{proof}

%
%

\section{Survival}
\label{S:surv}

In this section we prove Proposition~\ref{P:lambda_cCP(1)<=2}. Let $X=(X_t)_{t\geq0}$ be a cCP and let $\delta_0\in\mathcal{T}_\text{fin}$ be the configuration that equals 1 only at the origin. We say that $X$ \emph{survives} if 
\begin{align*}
\mathbb{P}^{\delta_0}\big[\exists t\geq0:X_t=\underline{0}\big]<1.
\end{align*}
The following result is known to hold for several processes. It is stated as \cite[Lemma~1]{SS08a} for an important class of cancellative processes. However, the cCP does not fit into this class and the definition of survival in the cited paper slightly differs from the one we are using here, so we provide a short proof below. Recall  that $\dot{\nu}$ is an invariant law of the cCP($\lambda,\delta$) that is defined as the long-time limit law of the process started in a non-trivial shift-invariant distribution, which is known to exist for $\lambda+\delta>0$ by \cite[Theorem~1.2\ \&\ Theorem~1.3]{BDD91}.
\begin{proposition}[Survival of the cCP]
\label{P:surv}
One has $\dot{\nu}\neq\delta_{\underline{0}}$ if and only if the \emph{cCP} survives.
\end{proposition}
\begin{proof}
We prove this statement using $\psib_\text{canc}$, the (multiplicative representation of the) cancellative duality function defined in (\ref{addANDcanc}). It is well-known that $\boldsymbol{\psi}_\text{canc}$ is informative, a fact that also follows from Proposition~\ref{P:genINFO}. Let $X=(X_t)_{t\geq0}$ be a cCP($\lambda,\delta$) ($\lambda,\delta\geq0,\ \lambda+\delta>0$) and let $x\in\mathcal{T}_\text{fin}$. If $\lambda,\delta>0$, then \cite[Theorem~1.2]{BDD91} implies that
\begin{align}
\label{BDD1}
\dot{\nu}\big(\{y:\vert x\wedge y\vert\ \text{is odd}\}\big)=\frac{1}{2}\mathbb{P}^x\big[X_t\neq\underline{0}\ \forall t\geq 0\big].
\end{align}
By the definition of $\boldsymbol{\psi}_\text{canc}$, (\ref{BDD1}) is equivalent to
\begin{align}
\label{theta1charac}
\int\boldsymbol{\psi}_\text{canc}(x,y)\ \di\dot{\nu}(y)=\mathbb{P}^x\big[\exists t\geq0:X_t=\underline{0}\big].
\end{align}
Choosing $x=\delta_0$ implies that $\dot{\nu}\neq\delta_{\underline{0}}$ if $X$ survives. On the other hand, if $X$ does not survive and $Y$ is a random variable with law $\dot{\nu}$, then (\ref{BDD1}) with $x=\delta_0$ implies that $\P[Y(0)=0]=1$ and the shift-invariance of $\dot{\nu}$ implies that $\P[Y(j)=0]=1$ for all $j\in\Z$. Hence $\dot{\nu}=\delta_{\underline{0}}$ as measures on $\mathcal{U}$ are characterised by their final dimensional marginals.
%
%

To complete the proof we consider the two special cases $\lambda=0$ and $\delta=0$. If $\lambda=0$, then $\delta>0$ and clearly $X$ does not survive while $\dot{\nu}=\delta_{\underline{0}}$. If $\delta=0$, then $\lambda>0$ and $X$ survives (one even has $\mathbb{P}^{\delta_0}[\exists t\geq0:X_t=\underline{0}]=0$) and $\dot{\nu}\neq\delta_{\underline{0}}$ by \cite[Theorem~1.3]{BDD91}.
\end{proof}

By Proposition~\ref{P:surv}, to prove Proposition~\ref{P:lambda_cCP(1)<=2}, it suffices to show that the cCP($\lambda,\delta$) does not survive when $\lambda\leq2\delta$. Let now $d=1$. Following \cite{Sud98} (compare the definition of $L$ in his Section~2), the idea for the proof of Proposition~\ref{P:lambda_cCP(1)<=2} is to construct a supermartingale applying Dynkin's formula to the function $g:\mathcal{T}_\text{fin}\setminus\{\underline{0}\}\to\mathbb{N}_0$ defined as
\begin{align}
\label{f:maxdist}
g(x):=\max\{i\in\mathbb{Z}:x(i)=1\}-\min\{i\in\mathbb{Z}:x(i)=1\}\qquad(x\in\mathcal{T}_\text{fin}).
\end{align}
In order to be able to apply Dynkin's formula one can ``reduce" the cCP to a finite state space similarly as in \cite[Proof of Lemma~3]{SS08b}. A full proof including the technical details is given below.

\begin{proof}[Proof of Proposition~\ref{P:lambda_cCP(1)<=2}]
Let $d=1$ and assume that $X$ is a cCP($\lambda,\delta$) with $\lambda\leq2\delta$. Using the $g$ from (\ref{f:maxdist}) we define $f:\mathcal{T}_\text{fin}\to\mathbb{N}_0$ as
\begin{align*}
f(x)=\begin{cases}
g(x)+4&\text{if}\ x\neq\underline{0},\\
0&\text{else},
\end{cases}
\qquad(x\in\mathcal{T}_\text{fin}).
\end{align*}
One then has that $G_\oplus f(x)\leq0$ for all $x\in\mathcal{T}_\text{fin}$, where $G_\oplus$ denotes the generator of the cCP from (\ref{gener}). To see this we first look at $x_{101},x_{11}\in\mathcal{T}_\text{fin}$ defined as
\begin{align*}
x_{101}(i)=\begin{cases}
1&\text{if}\ i\in\{0,2\},\\
0&\text{else},
\end{cases}\qquad
x_{11}(i)=\begin{cases}
1&\text{if}\ i\in\{0,1\},\\
0&\text{else},
\end{cases}
\qquad(x\in\mathbb{Z}).
\end{align*}
In the configuration $x_{101}$ the one at the origin reproduces with rate $\lambda$ to the left, increasing the function $f$ by one and it dies with rate $\delta$, decreasing $f$ by two. A reproduction to the right has no effect on $f$. By symmetry, an analogous statement holds for the one at $2\in\Z$ so that $G_\oplus f(x_{101})=2\lambda-4\delta$. For $x_{11}$ on the other hand, a reproduction of the one at the origin to the right reduces $f$ by one and its death reduces $f$ by only one, while a reproduction to the left again increases $f$ by one. Hence $G_\oplus f(x_{11})=-2\delta$. Let now $x\in\mathcal{T}_\text{fin}$ be an arbitrary configuration with at least two ones. As $f$ is shift-invariant, i.e.\ $f=f\circ\theta_i^{-1}$ for all $i\in\mathbb{Z}$, one has that $G_\oplus f(x)\leq G_\oplus f(x_{101})$ if $x$ has the form $010\ldots010$, $G_\oplus f(x)=G_\oplus f(x_{11})$ if $x$ has the form $011\ldots110$ and $G_\oplus f(x)\leq(G_\oplus f(x_{11})+G_\oplus f(x_{101}))/2$ if $x$ has the form $010\ldots110$ or $011\ldots010$. Note we had to use inequalities above as a death event of a one at the edge of a configuration reduces $f$ by the number of zeros ``to the inside" of this one, hence by at least two if there is a zero directly to the inside of the one. Finally we consider the special case $x=\delta_0$, in which with rate $2\lambda$ the lone individual reproduces (either to the left or to the right) and with rate $\delta$ it dies. Hence $G_\oplus f(\delta_0)=G_\oplus f(x_{101})=2\lambda-4\delta$, which was the reason to add the 4 in the definition of $f$. This completes the argument that $\lambda\leq2\delta$ implies that $G_\oplus f(x)\leq0$ for all $x\in\mathcal{T}_\text{fin}$.

The rest of the proof is a standard argument from the theory of continuous-time Markov chains, but, for the sake of completeness, we state it completely. Let $N\in\mathbb{N}$ be arbitrary and set $\tau_N:=\inf\{t\geq0:f(X_t)\geq N+4\}$. We claim that $M^N=(M_t^N)_{t\geq0}$ defined as
\begin{align*}
M_t^N:=f(X_{t\wedge\tau_N})-\int_0^{t\wedge\tau_N}G_\oplus f(X_s)\ \text{d}s\qquad(t\geq0)
\end{align*}
is a martingale. Let
\begin{align*}
\mathcal{T}_N&:=\{x\in\mathcal{T}_\text{fin}:x(i)=0\ \text{if}\ i\notin\{0,\ldots,N-1\}\}\cup\{x_N\},
\intertext{where}
x_N(i)&:=\begin{cases}
1&\text{if}\ i\in\{0,N\},\\
0&\text{else},
\end{cases}
\qquad(i\in\mathbb{Z}).
\end{align*}
By shifting every $x\in\mathcal{T}_\text{fin}$ so that its leftmost 1 lies at the origin we can construct a continuous-time Markov chain $Y=(Y_t)_{t\geq0}$ on the finite state space $\mathcal{T}_N$ so that
\begin{align*}
M_t^N=f(Y_t)-\int_0^t G_\oplus f(Y_s)\ \text{d}s\qquad(t\geq0).
\end{align*}
As a continuous-time Markov chain on a finite state space $Y$ is a Feller process and Dynkin's formula implies that $M^N$ is indeed a martingale.

As $G_\oplus f(x)\leq0$ for all $x\in\mathcal{T}_\text{fin}$ we conclude that $M^s=(f(X_{t\wedge\tau_N}))_{t\geq0}$ is a uniformly integrable supermartingale and the martingale convergence theorem implies that $M^s$ converges almost surely and in $L_1$ to a random variable $M_\infty$. The random variable $M_\infty$ is supported on $\{0,N+4\}$ as $M_\infty\in\{1,\ldots,N+3\}$ would imply that there exists a $t_0\geq0$ such that $M_t^s=M_{t_0}^s\in\{1,\ldots,N+3\}$ for all $t\geq t_0$, which has probability zero. Hence
\begin{align*}
4=\mathbb{E}^{\delta_0}[f(X_0)]\geq\mathbb{E}[M_\infty]=(N+4)(1-\mathbb{P}(M_\infty=0))
\end{align*}
and we conclude that
\begin{align*}
\mathbb{P}^{\delta_0}(\exists t\geq0:X_t=\underline{0})\geq \mathbb{P}^{\delta_0}(\exists t\leq\tau_N:X_t=\underline{0})=\mathbb{P}(M_\infty=0)\geq\frac{N}{N+4}.
\end{align*}
As $N$ was arbitrary it follows that $\mathbb{P}^{\delta_0}(\exists t\geq0:X_t=\underline{0})=1$ and Proposition~\ref{P:surv} implies that $\dot{\nu}=\delta_{\underline{0}}$. This establishes that $\lambda_\text{cCP}\geq2$.
\end{proof}

\appendix

\section{Parity indeterminacy}
\label{S:App}

In this appendix we restate \cite[Corollary~9]{SS08a} in a more accessible form and show how it can be derived from the somewhat less accessible formulation in  \cite{SS08a}. Then we show how this result implies Lemma~\ref{L:indeter}. 
Recall from Section~\ref{S:PoI} and Section~\ref{S:basicDUAL} the definitions of the operator $\oplus$ (addition modulo 2), of $\mathcal{T}$, the space all functions from $\Z^d$ to $T=\{0,1\}$, of $\mathcal{T}_\text{fin}\subset\mathcal{T}$, and of the cancellative duality function $\psib_2$. Let $\Ai$ be the set of all matrices of the form $A=(A(i,j))_{i,j\in\Z^d}$ with $A(i,j)\in\{0,1\}$ for all $i,j\in\Z^d$ and $\sum_{i,j}A(i,j)<\infty$. For $A\in\Ai$ and $x\in\mathcal{T}$, we define $Ax\in\mathcal{T}_\text{fin}$, corresponding to the usual matrix-vector multiplication, as
\begin{align*}
Ax(i):=\bigoplus_{j\in\Z^d}\big(A(i,j)\cdot x(j)\big)\qquad(i\in\Z^d),
\end{align*}
where $\cdot$ denotes the usual product in $\R$. Let $A^\dgg(i,j):=A(j,i)$ denote the adjoint of $A$. We will be interested in an interacting particle system $X=(X_t)_{t\geq 0}$ with state space $\mathcal{T}$, that jumps from its current state $x$ as
\be\label{xAx}
x\mapsto x\oplus Ax\quad\mbox{with rate}\quad a(A),
\ee
where $(a(A))_{A\in\Ai}$ are non-negative rates and the operator $\oplus$ has to be interpreted in a pointwise sense, as well as the interacting particle system $Y=(Y_t)_{t\geq 0}$ that jumps as
\begin{align*}
y\mapsto y\oplus A^\dgg y\quad\mbox{with rate}\quad a(A).
\end{align*}
In order for these interacting particle systems to be well-defined, we assume that
\be\label{sumco}
\sup_{i\in\Z^d}\sum_{A\in\Ai}a(A)\big|\{j:A(j,i)=1\}\big|<\infty
\quand
\sup_{i\in\Z^d}\sum_{A\in\Ai}a(A)\big|\{j:A^\dgg(j,i)=1\}\big|<\infty.
\ee
Recall from Section~\ref{S:weak} that $\vert z\vert:=\vert\{i\in\Z^d:z(i)=1\}\vert$ ($z\in\mathcal{T}$). It is shown in \cite{SS08a} that under condition (\ref{sumco}), the processes $X$ and $Y$ are well-defined and satisfy the duality relation

\be\label{XYdual}
\P\big[|X_tY_0|\mbox{ is odd}\big]=\P\big[|X_0Y_t|\mbox{ is odd}\big]\qquad(t\geq 0)
\ee
whenever $X$ and $Y$ are independent and either $|X_0|$ or $|Y_0|$ is a.s.\ finite.

We will restate \cite[Corollary~9]{SS08a}, which gives sufficient conditions for the left-hand side of (\ref{XYdual}) to be close to $1/2$. We assume that the rates are translation invariant in the sense that
\be\label{trans}
a(\theta_iA)=a(A)\qquad(i\in\Z^d,\ A\in\Ai),
\ee
where $\theta_iA$ denotes the ``translated'' matrix $(\theta_iA)(j,k):=A(j-i,k-i)$ $(j,k\in\Z^d)$. By definition, we say that a state $x\in\mathcal{T}$ is \emph{$X$-nontrivial} if
\be\label{Xnontriv}
\P^x\big[\big(X_t(i)\big)_{i\in\De}=\big(z(i)\big)_{i\in\De}\big]>0
\quad\mbox{for all $t>0$, finite $\De\sub\Z^d$, and }\big(z(i)\big)_{i\in\De}\in\{0,1\}^\De.
\ee
We fix a finite subset $\Bi\sub\Ai$ such that $a(B)>0$ for all $B\in\Bi$ and we define, for $x\in\mathcal{T}$,
\begin{align*}
\|x\|_\Bi:=\big|\big\{i\in\Z^d:\exists y\in\mathcal{T}\mbox{ and }B\in\Bi\mbox{ s.t.\ }\psib_2\big(x,(\theta_iB)y\big)=1\big\}\big\vert.
\end{align*}
With these definitions, \cite[Corollary~9]{SS08a} can be restated as follows. Recall the definition of the (pointwise) minimum operator $\wedge$ from Section~\ref{S:weak}.

\bp[Parity indeterminacy]
Let\label{P:indeter} $X$ be started in a shift-invariant initial law that is concentrated on $X$-nontrivial configurations. Then for each $\eps>0$ and $t>0$, there exists an $N<\infty$ such that
\be\label{indeter}
\left|\P\big[|X_t\wedge y|\mbox{ is odd}\,\big]-\frac{1}{2}\right|\leq\eps
\ee
for all $y\in\mathcal{T}_\text{\emph{fin}}$ with $\|y\|_\Bi\geq N$.
\ep

\begin{proof}
This is a simple reformulation of \cite[Corollary~9]{SS08a}. There, it is proved that if $y_n\in\mathcal{T}_\text{fin}$ satisfy $\|y_n\|_\Bi\to\infty$, then $\P\big[|X_t\wedge y_n|\mbox{ is odd}\big]\to\ha$. To see that this implies the claim of Proposition~\ref{P:indeter}, note that if the claim would be false, then there exists an $\eps>0$ such that for all $n\geq 1$ one can find $y_n\in\mathcal{T}_\text{fin}$ with $\|y_n\|_\Bi\geq n$ such that the left-hand side of (\ref{indeter}) is $>\eps$, contradicting \cite[Corollary~9]{SS08a}.
\end{proof}
Applying Proposition~\ref{P:indeter} to the cancellative contact process we obtain Lemma~\ref{L:indeter}.
\begin{proof}[Proof of Lemma~\ref{L:indeter}]
We first show that the jump rates of the cancellative contact process can be cast in the form (\ref{xAx}). Let $e_1,\ldots,e_d\in\Z^d$ denote the unit vectors and let $0\in\Z^d$ denote the origin. For $1\leq k\leq d$, we define $I^\pm_k\in\Ai$ by $I^\pm_k(i,j):=1$ if $(i,j)=(\pm e_k,0)$ and $I^\pm_{k}(i,j):=0$ otherwise. Also, we define $D\in\Ai$ by $D(i,j):=1$ if $(i,j)=(0,0)$ and $D(i,j):=0$ otherwise. Finally, we define rates $\big(a(A)\big)_{A\in\Ai}$ by
\begin{align*}
a(\theta_iI^\pm_{k}):=\la\quand a(\theta_iD):=\de\qquad(i\in\Z^d,\ 1\leq k\leq d),
\end{align*}
and $a(A):=0$ in all other cases. Clearly, these rates are translation invariant in the sense of (\ref{trans}) and satisfy the summability condition (\ref{sumco}). Also, a jump of the form $x\mapsto x\oplus(\theta_{-i}I^\pm_{k})x$ corresponds to a jump of the form $x\mapsto{\tt inf}^\oplus_{i,i\pm e_k}(x)$ in the notation of Section~\ref{S:PoI} and a jump of the form $x\mapsto x\oplus(\theta_{-i}D)x$ corresponds to a jump of the form $x\mapsto{\tt dth}_i(x)$, so the process defined by these rates is a cCP$(\la,\de)$. The claim of Lemma~\ref{L:indeter} will now follow from Proposition~\ref{P:indeter} provided we show that: (i)\ each configuration $x\neq\un 0$ is $X$-nontrivial and: (ii)\ we can choose $\Bi$ such that $\|y\|_\Bi=|y|$.

We start by proving (ii). We set $\Bi:=\{I_1^+\}$, where $I_1^+$ as defined above is one of the matrices corresponding to an infection next to the origin. Then $a(I_1^+)=\lambda>0$. Moreover,
\begin{align*}
\psib_2\big((\theta_{-i}I_1^+)x,y\big)=x(i)\cdot y(i+e_1)
\end{align*} 
and hence
\begin{align*}
y(i)=1\quad\text{if and only if}\quad e_1-i\in\{i\in\Z^d:\exists x\in\mathcal{T}\mbox{ and }B\in\Bi\mbox{ s.t.\ }\psib_2\big((\theta_iB)x,y\big)=1\big\},
\end{align*}
which shows that $\|y\|_\Bi=|y|$.

It remains to prove (i). Fix $x\in\mathcal{T}\setminus\{\un 0\}$, a finite set $\De\sub\Z^d$, and $\big(z(i)\big)_{i\in\De}\in\{0,1\}^\De$. Using the fact that $x\neq\un 0$ and $\la>0$, in a finite number of infection steps, we can infect each site in $\De\cup\{i\in\mathbb{Z}^d:\exists j\in\De:j\sim i\}$. Starting with the sites in $\De$ with the highest graph distance to $\mathbb{Z}^d\setminus\De$, we then can remove the infection from all sites $i$ such that $z(i)=0$ only using further infections, proving that the probability in (\ref{Xnontriv}) is positive for each $t>0$.
\end{proof}

The true strength of Proposition~\ref{P:indeter} lies in the fact that it can be applied even in situations where the definitions of $X$-nontriviality and the norm $\|y\|_\Bi$ are more complicated. In particular, \cite[Theorem~3]{SS08a} is based on an application of Proposition~\ref{P:indeter} in a situation where the $X$-nontrivial configurations are all $x\neq\un 0,\un 1$, and $\|y\|_\Bi=\big|\{(i,j):|i-j|=1,\ y(i)\neq y(j)\}\big|$.

\end{document}